\documentclass[reqno,12pt]{article}
\usepackage{a4wide}
\usepackage{amsmath,amscd} 
\usepackage{amssymb}
\usepackage{amsthm}
\usepackage[utf8]{inputenc} 

\usepackage{graphicx} 
\usepackage[all,cmtip]{xy} 
\usepackage{diagxy}
\usepackage[english, russian]{babel}

\numberwithin{equation}{section}

\newtheorem{theo}{Theorem}

\newtheorem{coro}{Corollary}
\newtheorem{prop}{Proposition}

\newtheorem{defi}{Definition}

\theoremstyle{remark}
\newtheorem{Remark}{Remark}

\def\al{\alpha}

\def\ze{\zeta}

\def\({\left(}
\def\){\right)}
\def\[{\left[}
\def\]{\right]}

\def\dd{\textup{d}}

\def\Qb{\overline{\mathbb Q}}

\newcommand{\N}{\mathbb{N}}
\newcommand{\Z}{\mathbb{Z}}
\newcommand{\Q}{\mathbb{Q}}
\newcommand{\R}{\mathbb{R}}
\renewcommand{\C}{\mathbb{C}}
\newcommand{\K}{\mathbb{K}}
\newcommand{\E}{{\bf E}}
\newcommand{\Qbar}{\overline{\mathbb Q}}
\newcommand{\Rplusetoile}{\mathbb{R}_+\etoile}
\newcommand{\Qplusetoile}{\mathbb{Q}_+\etoile}
\newcommand{\etoile}{^\star}
\newcommand{\Card}{{\rm Card}}

\newcommand{\eps}{\varepsilon}

\newcommand{\cinu}{\binom{\nu}{i}}
\newcommand{\ckj}{\binom{j}{k}}

\newcommand{\ckjmu}{\binom{j-1}{k}}

\newcommand{\Gh}{\widehat{\Gamma}}
\newcommand{\Ezero}{E_0}
\newcommand{\Erho}{E_\rho}
\newcommand{\Einf}{E_\infty}
\newcommand{\Nilsrho}{{\cal N}_\rho}
\newcommand{\Nilszero}{{\cal N}_0}
\newcommand{\Nilsinf}{{\cal N}_\infty}
\newcommand{\mrho}{m_\rho}

\newcommand{\calFbar}{\overline{\cal F}}
\newcommand{\calF}{{\cal F}}

\newcommand{\Lseul}{{\mathcal R}}	
\newcommand{\Lrho}{\Lseul_\rho}

\newcommand{\Lzero}{\Lseul_0}
\newcommand{\Linf}{\Lseul_\infty}
\newcommand{\Vinf}{{\mathcal L}_\infty}
\newcommand{\Vrho}{{\mathcal L}_\rho}
\newcommand{\Vloc}{{\mathcal L}_{\rm loc}}

\newcommand{\Invlapth}{\Vloc^{-1}}
\newcommand{\Invlapinf}{\Vinf^{-1}}
\newcommand{\Invlaprho}{\Vrho^{-1}}
\newcommand{\enhaut}{^\infty}
\newcommand{\asy}{{\mathcal A}_{\theta}}
\newcommand{\asymu}{{\mathcal A}_{\theta}^{-1}}
\newcommand{\asyeps}{{\mathcal A}\enhaut_{\theta+\eps}}
\newcommand{\asyepsmu}{{\mathcal A}^{\infty \, -1}_{\theta+\eps}}

\newcommand{\asymoinsepsmu}{{\mathcal A}^{\infty \, -1}_{\theta-\eps}}
\newcommand{\asypr}{{\mathcal A}\enhaut_{\theta'}}
\newcommand{\asyprmu}{{\mathcal A}^{\infty \, -1}_{\theta'}}
\newcommand{\tay}{{\mathcal A}^{(0)}}
\newcommand{\pro}{\kappa}
\newcommand{\Cseul}{{\mathcal T}}
\newcommand{\Czero}{\Cseul_0}
\newcommand{\Cinf}{\Cseul_\infty}
\newcommand{\matpro}{M_\kappa}
\newcommand{\matasy}{M_{\mathcal A^\infty_\theta}}
\newcommand{\matCzero}{  M _{\Czero}}
\newcommand{\matCinf}{  M _{\Cinf}}
\newcommand{\cloc}{c_{\rm loc}}
\newcommand{\cinf}{c_{\infty}}
\newcommand{\matVloc }{M_{\Vloc}}
\newcommand{\matVinf }{M_{\Vinf}}

\newcommand{\calD}{{\cal {D}}}  
\newcommand{\calDbar}{\overline{\cal {D}}}  
\renewcommand{\Re}{\textup{Re}\,}  
\renewcommand{\Im}{\textup{Im}\,}

\newcommand{\nvsing}{\Sigma}
\newcommand{\cercle}{{\cal C}(0,R)}

\newcommand{\deltarho}{\Delta_\rho}
\newcommand{\pointrho}{z_\rho}

\newcommand{\Arho}{A_\rho}
\newcommand{\GG}{{\bf G}}

\newcommand{\NGA}{{\rm NGA}}
\newcommand{\hada}{\star}

\newcommand{\Irho}{I_\rho}

\newcommand{\cz}{\C[z]}
\newcommand{\czN}{\C[z]_{<N}}
\newcommand{\Orho}{{\mathcal O}_\rho}
\newcommand{\Orhoch}{\widehat{\mathcal O}_\rho}

\newcommand{\Sinf}{{\mathcal S}_\infty}
\newcommand{\Srho}{{\mathcal S}_\rho}
\newcommand{\Sbarzero}{\overline{{\mathcal S}}_0}
\newcommand{\Sbarth}{\overline{{\mathcal S}}_{[\theta]}}
\newcommand{\Sbarthrho}{\overline{{\mathcal S}}_{[\theta]}^\rho}
\newcommand{\Sbarinfrho}{\overline{{\mathcal S}}_{\infty}^\rho}

\newcommand{\Sbarinf}{\overline{{\mathcal S}}_{\infty}}

\newcommand{\Sbarrho}{\overline{{\mathcal S}}_{\infty,\rho}}

\newcommand{\Zneg}{\Z_{\leq 0}}

\newcommand\tra{ \ ^t  }

\newcommand{\czM}{\C[z]_{<M}}
\newcommand{\zeromu}{\{0,\ldots,\mu\}}
\newcommand{\Span}{{\rm Span}}
\newcommand{\rk}{{\rm rk}\, }
\newcommand{\cut}{\Omega_\theta}
\newcommand{\cutbar}{\overline{\Omega}_\theta}

\newcommand{\antie}{\foreignlanguage{russian}{\CYREREV}}

\newcommand{\calV}{\mathbf{V}}
\newcommand{\hol}{{\cal H}}

\newcommand{\etanv}{\ae}
\renewcommand{\ae}{\mathfrak{f}} 
\newcommand{\ensembleae}{{\bf \antie }\ }
\newcommand{\ensembleaesansesp}{{\bf \antie }}

\newcommand{\epsz}{{\eps_0}}

\newcommand{\Bbarinf}{\overline{\mathcal B}_\infty}
\newcommand{\Bbarz}{\overline{\mathcal B}_0}
\newcommand{\Sto}{A}
\newcommand{\alin}{a}
\newcommand{\calP}{{\mathcal{P}}}

\begin{document}

 \selectlanguage{english}

\title{Microsolutions of differential operators and values of arithmetic Gevrey series}
\date\today
\author{S. Fischler and T. Rivoal}

\maketitle

\begin{abstract} We continue our investigation of $E$-operators, in particular their connection 
with $G$-operators; these differential operators are fundamental in understanding 
the diophantine properties of Siegel's  $E$ and $G$-functions. 
We study in detail  microsolutions (in Kashiwara's sense) of Fuchsian differential operators, and apply this to the construction of basis of solutions at $0$ and $\infty$ of any  $E$-operator from microsolutions of a $G$-operator; this provides a constructive proof of a theorem of Andr\'e. 
We also focus on the arithmetic nature of connection constants 
and Stokes constants between different bases of solutions of $E$-operators. For this, we introduce and study in 
details an arithmetic (inverse) Laplace transform that enables one to get rid of transcendental numbers inherent to 
Andr\'e's original approach. As an application, we define a set of special values of arithmetic 
Gevrey series, and discuss its conjectural relation with the ring of exponential periods.
\end{abstract}

 \section{Introduction}

In this paper, we continue our investigation of the arithmetic properties of 
certain differential operators related to $E$ and $G$-functions. 
Throughout the paper we fix a complex embedding of $\Qbar$ and  let $\N = \{0,1,2,\ldots\}$. To begin with, let us recall the following definition, essentially due to Siegel.
\begin{defi} \label{def:gfunc}
A $G$-function $G$ is a formal power series $G(z)=\sum_{n=0}^{\infty} a_nz^n$ 
such that the coefficients $a_n$ are algebraic numbers and there exists $C>0$ such that:
\begin{enumerate}
\item[$(i)$] the maximum of the moduli of the conjugates of $a_n$ is $\leq C^{n+1}$ for any $n$.

\item[$(ii)$] there exists a sequence of rational integers $d_n$, with $\vert d_n \vert \leq C^{n+1}$, such that 
$d_na_m$ is an algebraic integer for all~$m\le n$.

\item[$(iii)$] $G(z)$ satisfies a homogeneous linear differential equation with 
coefficients in $\Qb(z)$.
\end{enumerate}
\end{defi}
 An $E$-function is defined similarly, as $\sum_{n=0}^{\infty} \frac{a_n}{n!} z^n$ with 
the same assumptions $(i)$, $(ii)$, $(iii)$ (with $E(z)$ instead of $G(z)$ there).

\medskip
 
A minimal differential equation satisfied by a given $G$-function is called a $G$-operator. Any $E$-function is solution of 
an $E$-operator (not necessarily minimal) obtained as the Fourier-Laplace transform of a $G$-operator. In~\cite{YA1}, Andr\'e set the foundations of the theory of $E$-operators; he proved in particular that $0$ and $\infty$ are the only possible singularities, 
with rational exponents, and that 0 is a regular one. He also constructed two special bases of solutions of any $E$-operator at $0$ and $\infty$ respectively, and obtained a fundamental duality between these bases. A basic approach for constructing such bases is to lift the basis of solutions of the underlying $G$-operator through the Fourier-Laplace transform. However, it is well-known that there may not be enough such solutions to built a basis of the $E$-operator. A crucial ingredient in Andr\'e's construction is the notion of microsolutions (of the $G$-operator here) in the sense of  Kashiwara \cite{Kashiwara}, which 
enables one  to lift enough ``solutions'' of the $G$-operator. 

Our first theorem below is aimed at a better understanding of microsolutions of Fuchsian operators, of which $G$-operators are a special case. We shall explain in detail after the statement of the theorem the implications to Andr\'e's $E$-operators theory.  For any $\rho\in\C$, let $\Orho$ denote the ring of germs of functions holomorphic at $\rho$. Given $\theta\in\R$, let $\Orhoch$ denote the inductive limit as $\eps\rightarrow 0$ of the space of functions holomorphic on 
\begin{equation} \label{eqdisquecoupe}
\{z\in\C, \, \, 0 < |z-\rho| < \eps, \,   \, -\theta-\pi < \arg (z-\rho)  < -\theta+\pi\}.
\end{equation}

Given also a differential operator $\calD \in \C[z, \frac{\dd}{\dd z}]$, we denote by $\Sigma $ the set of all finite singularities of $\calD$, and by $\cut$ the   cut plane   obtained from $\C$ by removing the union of all closed half-lines of direction $- \theta + \pi $ starting at elements of $\Sigma$.

\begin{theo} \label{thintroun}
Let  $\calD \in \C[z, \frac{\dd}{\dd z}]$ and $\theta\in\R$ be such that $\arg(\rho-\rho')  \not\equiv  - \theta    \bmod \pi$ whenever $\rho,  \rho' \in \Sigma$ are distinct. Assume that $\infty$ is a regular singularity of $\calD$.

For any $\rho\in\Sigma$, let $f_\rho\in\Orhoch$. Then there exists a function $f$, holomorphic on the cut plane $\cut$, such that:
\begin{itemize}
\item $\calD f$ is a polynomial.
\item For any $\rho\in\Sigma$, we have $f-f_\rho\in\Orho$.
\end{itemize}
Moreover if $f$ and $\tilde f$ have these properties then $f-\tilde f$ is a polynomial.
\end{theo}

We refer to \S 1.4 of \cite{Pham} for an analogous result.
We shall prove Theorem \ref{thintroun} in \S \ref{sec3},   with an explicit upper bound on the degree of $\calD f$. The interest of this result is twofold.  On the one hand, given $\rho \in \Sigma$ any microsolution at $\rho$  (see below) can be represented by a fonction holomorphic on $\cut $ and such that $\calD f$ is a polynomial; this is of special interest to us when $\calD$ is a $G$-operator, because $f$ can therefore be expressed in terms of $G$-functions using the Andr\'e-Chudnovski-Katz Theorem. On the other hand, microsolutions at all finite singularities of $\calD$ can be glued together to produce a global function $f$ (as in \cite{Pham}).

Theorem \ref{thintroun} will be used in this paper only when $\calD$ is Fuchsian, so that we did not  try to remove the assumption that $\infty$ is a regular singularity of $\calD$. 
The proof of Theorem \ref{thintroun} is based on a linear bijective map $\pro$ (provided by analytic continuation) from 
$$
\Sinf = 
\{f \mbox{ holomorphic on $\cut$, } \calD f \in\cz\}/\cz
$$
to $\oplus_{\rho\in\Sigma} \Srho$. Here 
$$
\Srho  =
\{f\in\Orhoch, \, \, \calD f \in \Orho\}/\Orho
$$
  is the space of microsolutions of $\calD$ at $\rho$ as considered by  Kashiwara \cite{Kashiwara}. When $\calD$ is a $G$-operator, we prove also that $\pro$ and $\pro^{-1}$ are represented (in algebraic  bases) by matrices of which all coefficients are values of $G$-functions at algebraic points.
  
  \bigskip
  
We apply Theorem \ref{thintroun} more specifically to the following setting. Let   $\calFbar : \C[z,\frac{\dd}{\dd z}] \rightarrow  \C[x,\frac{\dd}{\dd x}]$ denote the Fourier-Laplace transform of differential operators, i.e. the morphism of $\C$-algebras defined by $\calFbar(z) =   \frac{\dd}{\dd x}$ and $\calFbar(\frac{\dd}{\dd z}) = -x$. Let  $\calDbar$ be an $E$-operator, that is  $\calFbar \calD$ for some $G$-operator  $\calD$ (see \cite{YA1}). We denote by $\Sbarzero$, resp. $\Sbarinf$, the space of local solutions of $\calDbar$ at the origin, resp. of formal solutions at infinity. We recall \cite{YA1} that a {\em   Gevrey series of order $s\in\Q$ of arithmetical type} is a power series $\sum_{n=0}^\infty a_n z^n$ with $a_n \in \Qbar$ such that  $\sum_{n=0}^\infty \frac{ a_n }{n!^s} z^n$ is a $G$-function, except maybe that it doesn't satisfy the holonomy assumption. Provided it is holonomic, with $s=0$, resp. $s=-1$, resp. $s=1$, this is  a $G$-function, resp. an $E$-function, resp. an  \antie-function \cite{YA1}.
These series form a differential algebra denoted by $ \Qbar\{z\}^A_{s}$.  A finite sum 
$$\sum_{\al\in S}\sum_{j\in T} z^{\al} (\log(z))^j \sum_{k\in K}\lambda_{\al,j,k} h_{\al,j,k} (z)$$
with $S\subset\Q$, $T,K\subset\N$, $\lambda_{\al,j,k}\in\C$ and  $h_{\al,j,k} (z) \in   \Qbar \{z\}^A_s$, is called in \cite{YA1} a {\em  Nilsson-Gevrey series} of order $s $ of arithmetical type. This differential algebra is   denoted by $\NGA\{z\}_s$. We denote by  $\NGA\{z\}^{\Qbar}_s$ the subset obtained by restricting to algebraic coefficients $\lambda_{\al,j,k}$. 

Andr\'e proved \cite{YA1} that $\Sbarzero$ has a basis $\Bbarz$ in $\NGA\{x\}^{\Qbar}_{-1}$, i.e. consisting of Nilsson-Gevrey series of order $-1 $ of arithmetical type with algebraic coefficients $\lambda_{\al,j,k}$.  Denoting again by $\Sigma$ the set of finite singularities of $\calD$, he has proved also that $\Sbarinf$ has a basis $\Bbarinf$ in $\oplus_{\rho\in\Sigma}e^{\rho x} \NGA\{1/x\}^{\Qbar}_{ 1}$.  We refer to Theorem \ref{thYA} in \S \ref{subsecnvpreuve} below for the precise statement of Andr\'e's result, including duality relations between $\Bbarz$ and $\Bbarinf$.

In this paper we give a new proof of Andr\'e's result, in a more constructive way. Andr\'e's proof is based on Laplace transform, which connects  solutions of $\calDbar$ to microsolutions of $\calD$ and enables him to apply fundamental results on $G$-functions. We revisit his approach in several directions. First we factor the inverse Laplace transform as $\Cseul \circ \Lseul$, where $\Cseul$ is a transcendental part (involving values of the Gamma function and its derivatives at rational points), and $\Lseul$ is a rational part (with only rational coefficients in suitable bases); $\Lseul$  is strongly related to Manjra-Remmal's operator \cite{ManjraRemmal} but is not exactly the same. Moreover $\Cseul$ commutes with differential operators, whereas $\Lseul$ behaves like the inverse  Laplace transform with respect to differential operators. We prove that 
$\Lseul$ induces (for $\rho \in \Sigma$) an explicit linear isomorphism defined over $\Qbar$ between the space $\Srho$ of microsolutions of $\calD$ at $\rho$ and the space of solutions of $\calDbar$ in $e^{\rho x} \NGA\{1/x\}^{\Qbar}_{ 1}$. Since microsolutions of $\calD$ can be represented in terms of $G$-functions (see the remark following Theorem \ref{thintroun} above), this provides an effective construction of a basis $\Bbarinf$ of $\Sbarinf$ as in \cite{YA1}. We obtain also, in the same effective way, a basis $\Bbarz$ of $\Sbarzero$. A major  difference  with Andr\'e's paper is the use of $\Lseul$, which is a formal operator,  instead  of the Laplace transform itself and operational calculus to give it a meaning when it is divergent. This enables us to get rid completely of $\Gamma$ values, and to work only with algebraic coefficients $ \lambda_{\al,j,k}$; we also 
obtain a new proof of the formal version of Theorem 2.2 of \cite{Bresil}.
Another important difference is that our study includes the case of integer exponents of $\calD$, whereas Andr\'e uses a trick 
to avoid it (see p.~734 of \cite{YA1}).

\bigskip

Our point of view enables us to study from an arithmetic point of view  the bijective linear map $\asy: \Sbarzero \rightarrow \Sbarinf$ given by asymptotic expansion in a large sector bisected by $\theta$; here $\theta\in\R$ is a fixed non-anti-Stokes direction, and throughout the paper a large sector bisected by $\theta$ is a sector of the form $\theta-\frac{\pi}{2}-\eps < 
\arg(x) <  \theta+\frac{\pi}{2}+\eps$ with $ \eps>0$. We refer to \cite{Ramis} or \cite{ateo} for the definition of  asymptotic expansion; in   other words, $ \asymu$ is Ramis' 1-summation in the direction $\theta$ (which coincides in this case with Borel-Laplace summation). Let $\matasy$ denote the matrix of $\asy$ in bases $\Bbarz$ and $\Bbarinf$ as above. The main result of \cite{ateo} asserts that $\matasy$ has coefficients in the $\GG$-module ${ \bf S}$ generated by all values at rational points of $\Gamma$ and its derivatives, where $\GG$ is the ring of values at algebraic points of analytic continuations of $G$-functions (see \cite{gvalues}). We recall that ${\bf S}$ is also a ring. We shall prove in \S \ref{subsecarith} that  $\matasy$ has a non-zero algebraic determinant, so that $ \matasy^{-1}$ has also coefficients in ${\bf S}$. This enables us to obtain the following result: 
 
 \begin{theo} \label{thnvcor}
Let $\ae(1/x)$ be an \antie-function, and $\theta\in\R$ be a non-anti-Stokes direction. Then there exist $\Sto \geq 1$  and $F_{\alin}\in  \NGA\{ x\}^{\Qbar}_{ 1}$, $\varpi_{\alin}\in{\bf S}$ (for $1\leq \alin\leq\Sto$) such that
$$ \asymu \ae(1/x) = \sum_{\alin=1}^{\Sto}  \varpi_{\alin } F_{\alin }(x).$$
\end{theo}

With respect to \cite{YA1}, the new feature in this corollary  is that $\varpi_{\alin}\in{\bf S}$.  
To state our next result, we recall that a Stokes matrix (relative to an $E$-operator $\calDbar$) is the matrix of change of coordinates between $  \asymu \Bbarinf$ and $  \asyprmu \Bbarinf$, where $\theta$ and $\theta'$ are non-anti-Stokes directions and $\Bbarinf$ is a basis of $\Sbarinf$ as above. In general 
$\theta$ and $\theta'$ are assumed to be one anti-Stokes direction apart from each other, but we shall not need this assumption here. 

\begin{coro} \label{corintrode} Any Stokes matrix (relative to an $E$-operator $\calDbar$) has coefficients in ${\bf S}$ and a non-zero algebraic determinant.
\end{coro}

Finally, we define a very large class of values of special functions.

\begin{defi} We denote by $\calV$ the ${\bf S}$-module generated by $\E$.
\end{defi}

We recall \cite{ateo} that $\E$ is the set of values at algebraic points of $E$-functions; it is a ring. Therefore  $\calV$ is also a ring; it consists in all finite sums of products $G(1)E(1)\Gamma(r)\gamma^j$ where $G$ is the analytic continuation of a $G$-function, $E$ is an $E$-function, $r\in\Q\setminus\Z_{\leq 0}$, $j\in\N$, and $\gamma$ is Euler's constant (see \cite{gvalues} and \cite{ateo}). For instance $\calV$ contains all values at algebraic points of Airy's oscillating integral, and  that of Bessel's functions $J_\alpha(z)$ with $\alpha\in\Q$.
Actually it  contains all values at algebraic points of 
  Nilsson-Gevrey series  of any order   of arithmetical type with algebraic coefficients $\lambda_{\al,j,k}$ (up to 1-summation in any direction in the case of divergent series), and especially  values of 
all generalized hypergeometric series.

We shall deduce from the duality between $E$- and \antie-functions a dual expression of $\calV$. Let us denote by \ensembleae the set of all complex numbers $\ae_\theta(\xi)$ where $\xi\in\Qbar\etoile$  and $\ae$ is an \antie-function; here $\theta = \arg(\xi)$ and  $\ae_\theta =  \asymu \ae$ is Ramis' 1-summation of $\ae$ in the direction $\theta$ if $\theta$ is not anti-Stokes, and  $\ae_\theta = {\asyeps}^{-1}\ae$ for any small $\eps>0$ if $\theta$ is  anti-Stokes (this is independent from the choice of such an $\eps$). 
 Andr\'e has conjectured \cite{YAJTNB2003} that Siegel-Shidlovskii theorem on values of $E$-functions has an analogue for \antie-functions.
We shall prove in \S \ref{subsec43} that \ensembleae is a ring which contains   algebraic numbers,  Gompertz' s constant $\int_0 ^{+\infty}  \frac{e^{-t}}{1+t} \dd t$, and $\sqrt{\pi} {\rm Ai}(z)$ for any $z\in\Qbar$ where ${\rm Ai}(z)$ is Airy's oscillating integral. Using Theorem \ref{thnvcor} and the results of \cite{ateo} we obtain the following dual characterization of $\calV$.

\begin{theo} \label{thintrotr}
The ${\bf S}$-module generated by the numbers $e^\rho \chi$, with $\rho\in\Qbar$ and $\chi\in$ \ensembleaesansesp, is equal to $\calV$.
\end{theo}

\medskip

We believe that $\calV$ is related to the ring $\calP_e$ of exponential periods, defined by Kontsevich in the last section of \cite{KZ}. An exponential period is an absolutely convergent integral $\int_\Omega f(x) \exp(g(x)) \dd x$ where $n \geq 1$, $x = (x_1,\ldots,x_n)$, $\Omega \subset \R^n$ is a semi-algebraic domain (i.e., it is defined by polynomial inequalities and/or equalities with algebraic coefficients), and $f$, $g$ are algebraic functions. Restricting to $g=0$ yields the ring $\calP$ of periods, in the sense of Kontsevich-Zagier \cite{KZ}. In view of the Bombieri-Dwork conjecture, it is natural to ask whether $\GG = \calP[1/\pi]$ (see \S 2.2 of \cite{gvalues}). In the present setting, the corresponding question would be whether $\calV  = \calP_e[1/\pi]$ (see \S \ref{subsec43}).

\bigskip

The structure of this paper is as follows. In \S \ref{sec2}  we factor the (inverse) Laplace transform as announced above: we construct the operators $\Cseul$ and $\Lseul$, and study their properties. This part of the paper does not involve any specific class of differential operators; it would be interesting to find other situations where it could be used. In \S \ref{sec3} we study the microsolutions at all finite singularities of a differential operator $\calD\in\C[z, \frac{\dd}{\dd z}]$ of which $\infty$ is a regular singularity, and prove Theorem \ref{thintroun}. When $\calD$ is Fuchsian, we also explain how the analytic continuation map $\pro : \Sinf\stackrel{\sim}{\to}  \oplus_{\rho\in\Sigma} \Srho$ corresponds, through Laplace transform, to a duality between  the solutions of $\calFbar\calD$ at 0 and at $\infty$, given by asymptotic expansion. Using either $\pro$, of $\Lseul$ and $\Cseul$, we get two new constructions of Andr\'e's extension of Laplace transform (based on operational calculus). We relate them and   sum up our results in a commutative diagram (see \S \ref{secdiagLG}).  
At last we apply our results in \S \ref{secspecialvalues} to the case where $\calD$ is a $G$-operator. This enables us to obtain a new constructive proof of Andr\'e's duality theorem, and to prove the arithmetic results announced in this introduction. At last, we work out in \S \ref{subsec44} a simple example related to Gompertz' constant, to illustrate the whole situation.

\section{A rational version of the  inverse  Laplace transform} \label{sec2}

In this section we define formal operators $\Lrho$ for $\rho\in\C\cup\{\infty\}$  (see \S \ref{subsec11}), which behave like the inverse Laplace transform with respect to differential operators and involve only rational coefficients. Their properties are stated in \S \ref{subsec23}, and proved in \S \ref{subsec12}. Then we define in \S \ref{subsec25} formal operators $\Czero$ and $\Cinf$ which involve values of derivatives of $\Gamma$ and allow us to factor the Laplace transform and its inverse (see Propositions \ref{proplienlaplaceinverserho},   \ref{proplienlaplaceinverseinf}, and \ref{prop7}) and to write down explicitly the implicit formulas obtained by  Andr\'e using operational calculus (see Remark \ref{rem2} at the end of \S \ref{subsec25}).

\subsection{Notation} \label{subsec10}

Until the end of \S \ref{subsec12}, we consider   formal variables $z$ and $x$.

Let $\rho\in\C$. We denote by $\Erho$ the set of all (formal) functions that can be written as
\begin{equation}\label{eqdeferho}
f(z-\rho) = \sum_{\al\in S} \sum_{j\in T} \sum_{n=0}^{\infty} c_{\al,j,n} (z-\rho)^{n+\al}(\log(z-\rho))^j,
\end{equation}
where $S\subset \C$ and $T\subset\N$ are finite subsets, and $c_{\al,j,n} \in\C$ for any $\al$, $j$, $n$.
If we assume for a given $f$ that $S$ and $T$ have the least possible cardinality (so that  $\al - \al' \not\in\Z$ for any distinct elements $\al,\al'\in S$) and that  for any $\al \in S$ there exists $j \in T$   such that $c_{\al,j,0}\neq 0$, then the expansion \eqref{eqdeferho} is
uniquely determined by $f$. We say that $f$ has   {\em  coefficients  and exponents in a subfield} $\K$ of $\C$ if all $c_{\al,j,n} $ and all $\al$ involved in this expansion  belong to $\K$.

Of course $f(z-\rho)\in\Erho$ if, and only if, $f(z) \in \Ezero$.

   If all power series $\sum_{n=0}^{\infty} c_{\al,j,n} (z-\rho)^{n }$ have   positive radii of convergence, the function \eqref{eqdeferho} is said to belong to the Nilsson   class; we denote by $\Nilsrho$ this subset of $\Erho$. We have $\NGA\{z-\rho\}_s \subset \Erho$ for any $s\in\Q$, and $\NGA\{z-\rho\}_s \subset \Nilsrho$ if $s\leq 0$.

\bigskip

  In the same way, we also let $\Einf$ denote the set of all formal functions
\begin{equation}\label{eqdefeinf} 
f(x) = \sum_{\al\in S} \sum_{j\in T} \sum_{n=0}^{\infty} c_{\al,j,n}x^{-n-\al-1}(\log(1/x))^j,
\end{equation}
with the same remarks; $\Nilsinf$ is the corresponding Nilsson class.

\subsection{Definition of the operators $\Lrho$} \label{subsec11}

To begin with, let us define a {\em rational  inverse  Laplace transform} $\Lrho$ as follows, for $\rho\in\C\cup\{\infty\}$. For any $\al\in\C$ and any $j\in\N$ we first define
\begin{equation}\label{eqdeflseul}
\Lseul \Big( z^\alpha (\log(z))^j\Big) =  \sum_{k=0}^j \ckj \frac{\dd^{j-k}}{\dd y^{j-k}}\Big(\frac{\Gamma(1-\{y\})}{\Gamma(-y)}\Big)_{|y=\al} x^{- \al-1}(\log(1/x))^k.
\end{equation}
Here we let $\{u+iv\} = {\rm Frac}(u)+iv$ for $u,v\in\R$, where ${\rm Frac}(u) \in[0,1)$ is the fractional part of $u$; moreover 
the $(j-k)$-th derivative of $\frac{\Gamma(1-\{y\})}{\Gamma(-y)}$ is computed at the right of $\al$ (i.e., as $y \rightarrow \al$, $\Re y > \Re \al$) if $  \Re \al\in\Z$ (otherwise it is simply computed at $\al$).  

For any $\rho\in\C$ we define a  linear map $\Lrho : \Erho \rightarrow e^{\rho x} \Einf$ by letting
$$\Lrho\Big( \sum_{\al\in S} \sum_{j\in T} \sum_{n=0}^{\infty} c_{\al,j,n} (z-\rho)^{n+\al}(\log(z-\rho))^j\Big)
 =  \sum_{\al\in S} \sum_{j\in T} \sum_{n=0}^{\infty} c_{\al,j,n} e^{\rho x} \Lseul \Big( z^{n+\alpha} (\log(z))^j\Big)
  $$
with the same notation as in Eq. \eqref{eqdeferho}. These linear maps $\Lrho$, depending on $\rho\in\C$, are related to one another since
\begin{equation}\label{eqlienrho}
\Lrho (f(z-\rho)) = e^{\rho x} \Lzero(f(z))
\end{equation}
for any $f(z) \in \Ezero$.
Using Leibniz' rule  they can be written in a more compact way, upon noticing that 
\begin{equation}\label{eqastuce}
\Lrho\Big((z-\rho)^{ \al}(\log(z-\rho))^j\Big) =e^{\rho x}   \frac{\partial^{j }}{\partial y^{j }} \Big( \frac{\Gamma(1-\{y\})}{\Gamma(-y)}x^{- y -1}\Big)_{|y=\al} .
\end{equation}
We also define a linear map $\Linf : \Einf \rightarrow \Ezero$ by
$$\Linf\Big( \sum_{\al\in S} \sum_{j\in T} \sum_{n=0}^{\infty} c_{\al,j,n} z^{-n-\al-1}(\log(1/z ))^j\Big)
 =  \sum_{\al\in S} \sum_{j\in T} \sum_{n=0}^{\infty} c_{\al,j,n} (-1)^j  \Lseul \Big( z^{-n-\al-1}(\log(z)  )^j\Big).$$
This means that we let (formally) $\log(1/z) = - \log(z)$; with this convention $\Linf$ and  $\Lzero$ coincide on $z^\alpha \C[\log (z)]$ for any $\alpha\in\C$.    However $\Linf$ and $\Lzero$ are defined on different sets, namely $\Einf$ and $\Ezero$ respectively. 

\bigskip

\begin{Remark} \label{rem1}
Very similar formulas appear in \S 5 of \cite{ManjraRemmal}, where Manjra and Remmal define a formal Laplace transform with properties analogous to Propositions \ref{proplap} and  \ref{proprat} below. However they restrict to $\alpha\not\in\Z$, whereas the case of integer exponents is very important in our approach, and their version of the  Laplace transform depends on the choice of a certain square matrix $\Lambda$.
\end{Remark}

\subsection{Statement of the Properties of $\Lrho$} \label{subsec23}

The following propositions will be proved in \S \ref{subsec12} below. We recall that  $\calFbar : \C[z,\frac{\dd}{\dd z}] \rightarrow \C[x,\frac{\dd}{\dd x}]$ is  the Fourier-Laplace  transform of differential operators, i.e. the morphism of $\C$-algebras defined by $\calFbar(z) =   \frac{\dd}{\dd x}$ and $\calFbar(\frac{\dd}{\dd z}) = -x$.

\begin{prop}  \label{proplap}
For any $\rho\in\C\cup\{\infty\}$ and any $\calD\in \C[z,\frac{\dd}{\dd z}] $ we have
$$ \Lrho \circ \calD =  (\calFbar \calD) \circ \Lrho .$$
In particular, for any $f(z-\rho)\in\Erho$ we have
$$ \Lrho (z\, f(z-\rho)) =  \frac{\dd}{\dd x} \Lrho (f(z-\rho))  
\mbox{ and }  \Lrho \Big( \frac{\dd f}{\dd z} (z-\rho)\Big) = -x \Lrho (f(z-\rho)),$$
where $z-\rho$ should be understood as $1/z$ if $\rho=\infty$.
\end{prop}

These relations are satisfied by the usual inverse Laplace transform; they are the reason why $\Lrho$ is   a modified  inverse  Laplace transform. On the contrary,  the following rationality  property holds for $\Lrho$ but not for the usual  inverse  Laplace transform, and it is crucial in our approach.

\begin{prop} \label{proprat}
For any $\rho\in\C$, if $f(z-\rho) \in \Erho$  has coefficients and exponents in a subfield $\K$ of $\C$ then $e^{-\rho x} \Lrho (f(z-\rho))\in \Einf$ has coefficients and exponents in the same subfield.

In the same way, if $f(1/z)\in\Einf$   has coefficients and exponents in $\K$  then so does $\Linf( f(1/z))\in\Ezero$.
\end{prop}

Recall from \S \ref{subsec10} that an element of $\Erho$ or $\Einf$ is said to have  coefficients and exponents in $\K$ if it can be written as \eqref{eqdeferho} or  \eqref{eqdefeinf} with $S\subset \K$ and all coefficients $c_{\al,j,n}$ in $\K$.  

In particular, if we restrict $\Lrho$ to the subspace of $\Erho$ consisting of the functions  \eqref{eqdeferho}  with $S \subset \Q$, then we obtain a linear map defined over the rationals; the analogous property holds also for $\Linf$.

\bigskip

It will be important for us   that $\Lrho$ induces a bijective linear map 
$
\Erho/\C[[z-\rho]]
\rightarrow e^{\rho x}\Einf$, and also 
$\Linf : 
\Einf/\cz \rightarrow \Ezero$. This is the meaning of the following proposition.

\begin{prop} \label{propker}
Let $\rho\in\C$. Then $\Lrho : \Erho \rightarrow e^{\rho x}\Einf$ is surjective, and its kernel   is exactly the space  $\C[[z-\rho]]$ of formal power series $\sum_{n=0}^\infty c_n (z-\rho)^n$.

On the other hand, $\Linf : \Einf\rightarrow\Ezero$ is also  surjective, and its kernel   is exactly the space of polynomials $\cz$.
\end{prop} 

\bigskip

Finally, we shall use the fact that $\Lrho$ and $\Linf$ map Nilsson-Gevrey  series  of order $s\in\Q$ of arithmetical type with algebraic coefficients  to Nilsson-Gevrey  series  of order $s \pm 1$ of arithmetical type with algebraic coefficients. The same property holds for the usual inverse Laplace transform (see \cite{YA1}), except for the algebraicity of  coefficients; for this reason our proof is more direct.

\begin{prop} \label{proparith} 
If $f(z-\rho)\in \NGA\{z-\rho\}^{\Qbar}_s$ with $\rho\in\C$ and $s\in\Q$,  then $e^{-\rho x} \Lrho( f(z-\rho) ) \in \NGA\{1/x\}^{\Qbar}_{s+1}$.

On the other hand, if $f(1/z) \in \NGA\{1/z\}^{\Qbar}_s$ then $\Linf( f(1/z)) \in \NGA\{x\}^{\Qbar}_{s-1}$.
\end{prop}

\subsection{Proof  of these properties} \label{subsec12}
 
In this section we prove Propositions \ref{proplap}, \ref{proprat}, \ref{propker},  and \ref{proparith}. To begin with, recall that Pochhammer's symbol is defined by $(y)_p = y(y+1)\ldots(y+p-1)$. Given $y\in \C\setminus\N$  we denote by $n\in\Z$ the integer part of its real part, so that $\{y\} = y-n$ and  
\begin{equation} \label{eqconcretun}
\frac{\Gamma(1-\{y\})}{\Gamma(-y)} = \frac{\Gamma( - y +n+1)}{\Gamma(-y)} = 
\left\{ \begin{array}{l}
(-y)_{n+1} \mbox{ if } n\geq 0\\ 
\\
\frac1{(-y+n+1)_{-n-1}}\mbox{ if } n \leq -1
\end{array}\right.
\end{equation}
is a rational function of $y$, with rational coefficients (as long as $n$ is fixed). To differentiate this function of $y$ and then take $y=\al$ (or $y\rightarrow \al$, $\Re y > \Re \al$, if $  \Re \al\in\Z$), we may assume that  $\Re y$ and $\Re \al$ have the same integer part. Then the values of all derivatives of this rational function   belong  to $\Q(\al)$: all coefficients in Eq. \eqref{eqdeflseul} belong to $\Q(\al)$, and Proposition  \ref{proprat} follows immediately.

\bigskip

To prove Proposition \ref{proplap} for $\rho\in\C$, it enough to check that 
 for any $f(z-\rho)\in\Erho$ we have
$$ \frac{\dd}{\dd x} \Lrho (f(z-\rho)) = \Lrho (z\, f(z-\rho))
\mbox{ and } -x \Lrho (f(z-\rho)) = \Lrho \Big( \frac{\dd f}{\dd z} (z-\rho)\Big).$$
Let us   begin with the case $\rho=0$. We may assume 
$f(z ) =    z ^{ \al}(\log(z) )^j$, and Eq. \eqref{eqastuce} yields
\begin{align*}
 \frac{\dd }{\dd x} \Lzero(z^\al (\log(z))^j)  &=    \frac{\partial^{j }}{\partial y^{j }} \Big(\frac{\partial }{\partial x}  \Big( \frac{\Gamma(1-\{y\})}{\Gamma(-y)}x^{- y -1}\Big)\Big)_{|y=\al } \\
&=  \frac{\partial^{j }}{\partial y^{j }} \Big(  \frac{\Gamma(1-\{y \})}{\Gamma(-y-1)}x^{- y -2} \Big)_{|y=\al } \\
&=  \Lzero(z^{\al+1} (\log(z))^j) 
\end{align*}
and
\begin{align*}
-x  \Lzero (  z^\al (\log(z))^j ) &=
-x  \frac{\partial^{j }}{\partial y^{j }} \Big( \frac{\Gamma(1-\{y\})}{\Gamma(-y+1)}x^{- y  } \cdot  \frac{-y}{x} \Big)_{|y=\al } \\
&= \al   \frac{\partial^{j }}{\partial y^{j }} \Big( \frac{\Gamma(1-\{y\})}{\Gamma(-y+1)}x^{- y }\Big)_{|y=\al}  + j   \frac{\partial^{j-1 }}{\partial y^{j-1 }} \Big( \frac{\Gamma(1-\{y\})}{\Gamma(-y+1)}x^{- y }\Big)_{|y=\al }  \\
&= \al    \Lzero (  z^{\al-1} (\log(z))^j ) + j   \Lzero (  z^{\al-1} (\log(z))^{j-1} )\\
&=  \Lzero\Big(   \frac{\dd }{\dd z}\Big(  z^\al (\log(z))^j\Big)\Big) .
\end{align*}
These computations prove also Proposition \ref{proplap} for $\Linf$. To deduce the same property for $\Lrho$ with 
 any $\rho\in\C$, we use Eq. \eqref{eqlienrho}:
\begin{align*}
\frac{\dd }{\dd x}\Big( \Lrho(f(z-\rho))\Big)   &= 
 e^{\rho x} \Big(  \frac{\dd }{\dd x} L_0(f(z)) + \rho L_0(f(z)) \Big)
\\
& =  e^{\rho x}   \Big(  L_0(z f(z)  + \rho f(z)) \Big) 
\\
&= 
    \Lrho(z f(z-\rho))\quad  \mbox{ since } z = (z-\rho)+\rho
 \end{align*}
and
$$ \Lzero\Big(   \frac{\dd }{\dd z}\Big(  f(z-\rho) \Big)\Big) =
 e^{\rho x} \Lzero \Big(   \frac{\dd f }{\dd z} \Big) 
= - x e^{\rho x} \Lzero (   f (z) ) 
 = -x \Lzero(f(z-\rho)).$$

\bigskip

Let us prove  Proposition \ref{propker} now.  
  Let $\al\in\C\setminus\N$ and $j\in\N$. Then Eq. \eqref{eqdeflseul} asserts that $  \Lseul (  z ^{ \al}(\log(z) )^j )$ is equal to $  x^{- \al-1}$ times a polynomial in $\log(1/x)$ of degree at most $j$. Moreover the coefficient of degree $j$ of this polynomial is $\frac{\Gamma(1-\{\al \})}{\Gamma(-\al )} \neq 0$ since $\al\not\in\N$, so that $\Lrho$ induces a bijective degree-preserving map $(z-\rho)^\al \C[\log (z-\rho)] \rightarrow e^{\rho x} x^{- \al-1} \C[ \log(1/x)]$.

Now let us move to  non-negative integer values of $\al$. For $\al \in\N$, the function $y \mapsto \frac{\Gamma(1-\{y\})}{(y-\al)\Gamma(-y)}$ is holomorphic at the right of $\al$ and does not vanish at $y=\al$, and Leibniz' formula shows that 
$$\frac{\dd^{j-k}}{\dd y^{j-k}}\Big(\frac{\Gamma(1-\{y\})}{\Gamma(-y)}\Big)_{|y=\al }  = (j-k) \frac{\dd^{j-k-1}}{\dd y^{j-k-1}}\Big(\frac{\Gamma(1-\{y\})}{(y-\al)\Gamma(-y)}\Big)_{|y=\al } $$
for any $0\leq k \leq j$. Therefore this quantity   is zero for $k=j$ and non-zero for $k=j-1$, and  Eq. \eqref{eqdeflseul} yields  for $\al\in\N$:
$$\Lrho\Big((z-\rho)^{ \al}(\log(z-\rho))^j\Big) = j e^{\rho x} \sum_{k=0}^{j-1} \ckjmu \frac{\dd^{j-k-1}}{\dd y^{j-k-1}}\Big(\frac{\Gamma(1-\{y\})}{(y-\al)\Gamma(-y)}\Big)_{|y=\al } x^{- \al-1}(\log(1/x))^k.$$
 This shows that  $\Lrho$ induces in this case a   map $(z-\rho)^\al \C[\log (z-\rho)] \rightarrow e^{\rho x} x^{- \al-1} \C[ \log(1/x)]$ which decreases the degree by 1; it vanishes exactly on constant multiples of  $(z-\rho)^\al $. This concludes the proof of Proposition \ref{propker}.

\bigskip

Let us prove  Proposition \ref{proparith} now.  
By linearity we may assume that $f(z-\rho) = (z-\rho)^{\al} (\log (z-\rho))^j  h (z-\rho)$ with $\al\in\Q$, $j\in\N$, and $h\in \Qbar \{z-\rho\}^A_s$. Then
\begin{equation} \label{eqcalclrho}
\Lrho( f(z-\rho)) = j! e^{\rho x} \sum_{k=0}^j (y_{\al, j-k} \hada h)(1/x) x^{-\al-1} \frac{ (\log(1/x))^k}{k!}
\end{equation}
where $\hada$ is Hadamard coefficientwise product of formal series in $z$, and 
$$y_{\al,i}(z) = \sum_{n=0}^\infty \frac{1}{i!} \frac{\dd^{i}}{\dd y^{i}}\Big(\frac{\Gamma(1-\{y\})}{\Gamma(-y-n)}\Big)_{|y=\al } z^n.$$
  Now Eq. \eqref{eqconcretun} yields
$$ \frac{\Gamma(1-\{y\})}{\Gamma(-y-n)} = (-1)^{[y]+n+1} (\{y\})_{[y]+n+1}$$
for any $n \geq -[y]$ so that $y_{\al,0}(z)\in \Q\{z\}_{1}^A$, because $ (-1)^{[\al]+1} \mbox{}_2F_0\bigg(\begin{array}{cc}1,\; \{\al\} \\ \mbox{-} \end{array}\bigg| -z\bigg) \in \Q\{z\}_{1}^A$  since $\alpha\in\Q$ (see \cite{Andre}, Chapter I, \S 4.4).  It is not difficult   to prove more generally that  $y_{\al,i}(z)\in \Q\{z\}_{1}^A$
for any $\alpha\in\Q$ and any $i\in\N$. The first part of   Proposition \ref{proparith} follows easily, using Eq. \eqref{eqcalclrho}. 

To prove the second part, we proceed in the same way. Letting 
$$\tilde y_{\al,i}(z) = \sum_{n=0}^\infty \frac{1}{i!} \frac{\dd^{i}}{\dd y^{i}}\Big(\frac{\Gamma(1-\{y\})}{\Gamma(-y+n)}\Big)_{|y=\al } z^{-n}$$
we have for $\alpha\in\Q$, $j\in\N$ and $h(z) \in \Qbar \{1/z\}^A_s$:
$$
\Linf( z^{-\alpha-1}(\log(1/z))^j h(z) ) =(-1)^j j!   \sum_{k=0}^j (-1)^k  \big(\tilde y_{-\al-1, j-k} \hada h\big)( x) x^{ \al } \frac{ (\log  x)^k}{k!}
$$
where 
 $\hada$ is Hadamard's product of formal series in $1/z$, so that $(\tilde y_{-\al-1, j-k} \hada h)(x)$ is a power series in $x$. Now  $\tilde y_{\al,i}(z)\in \Q\{1/z\}_{-1}^A$
for any $\alpha\in\Q$ and any $i\in\N$, and  Proposition \ref{proparith} follows.

\subsection{Factorization of the inverse Laplace transform} \label{subsec25}

In this section we define differential operators $\Czero$ and $\Cinf$, the coefficients of which involve values of derivatives of the $\Gamma$ function. We prove that $\Cinf\Lrho f$ and $\Czero \Linf f$ are respectively the asymptotic expansion at infinity and the generalized Taylor expansion at 0 of inverse Laplace transforms of $f$ (where integration is performed along suitable paths). The Laplace transform itself is given by  $(\Cinf \circ \Lrho )^{-1} $ or  $(\Czero\circ  \Linf)^{-1} $; this enables one to make Andr\'e's proof of his duality theorem more explicit (see the remark at the end of this section).

\bigskip

For simplicity we let $\Gh(s) = 1/\Gamma(s)$. As in \S \ref{subsec11}  we let $\{u+iv\} = {\rm Frac}(u)+iv$ for $u,v\in\R$, where ${\rm Frac}(u) \in[0,1)$ is the fractional part of $u$. For any $\rho, \al\in\C$ and any $\nu\in\N$ we let 
$$\Cseul\Big(e^{\rho x} x^{-\al-1}(\log(1/x ))^\nu\Big) =  e^{\rho x} x^{-\al-1}  \sum_{i=0}^\nu (-1)^{\nu-i} \cinu \Gh ^{(\nu-i)}(1-\{\alpha\}) (\log(1/x ))^i$$
where the right handside can also be written as
$$e^{\rho x} \Big(\frac{\dd}{\dd y}\Big)^{\nu}\Big(  \Gh  (1-\{ y \}) x^{-y-1}\Big)_{|y=\alpha}.$$
As in the rest of this paper, all derivatives involving $\{y\}$ are computed at the right  of $\alpha$. 
We extend $\Cseul$ to a linear map $\Cinf : e^{\rho x} \Einf    \rightarrow e^{\rho x} \Einf$ by letting 
$$\Cinf \Big( \sum_{\al\in S} \sum_{j\in T} \sum_{n=0}^{\infty} c_{\al,j,n} e^{\rho x}  x^{-n-\al-1}(\log(1/x ))^j\Big)
= \sum_{\al\in S} \sum_{j\in T} \sum_{n=0}^{\infty} c_{\al,j,n}\Cseul \Big( e^{\rho x}  x^{-n-\al-1}(\log(1/x ))^j\Big)$$
and also to a  linear map $\Czero : \Ezero    \rightarrow \Ezero$ in the same way:
$$\Czero \Big( \sum_{\al\in S} \sum_{j\in T} \sum_{n=0}^{\infty} c_{\al,j,n} x^{ n+\al }(\log(x) )^j\Big)
= \sum_{\al\in S} \sum_{j\in T} \sum_{n=0}^{\infty} c_{\al,j,n} (-1)^j \Cseul \Big(  x^{ n+\al }(\log(1/x ))^j\Big).$$
In other words, as in \S \ref{subsec11}   we agree that, formally,  $\log(x) = - \log(1/x)  $.

It is not difficult to prove that $\Czero \circ \frac{\dd}{\dd x} =  \frac{\dd}{\dd x} \circ\Czero$ and $\Czero \circ x = x  \circ\Czero$  (where $x$ denotes multiplication with $x$), so that $\Czero$ commutes with any $\calD \in \C[x,\frac{\dd}{\dd x}]$; the same property holds for $\Cinf$.

\bigskip

Before we  can  state the relation with the inverse Laplace transform, we have to define the paths of integration (see    \cite[pp. 183--192]{DiP} or \cite{ateo}). Given $\theta\in\R$ and $\rho\in\C$, we consider the cut plane defined by $z \neq \rho$ and $-\theta-\pi < \arg (z-\rho)  < -\theta+\pi$. In this cut plane we denote by $\Gamma_\rho$ the following path: a straight line from $\rho + e^{i(-\theta-\pi)}\infty$ to $\rho$ (on one bank of the cut), a circle of radius essentially zero around $\rho$ (with $\arg (z-\rho)  $ increasing from $-\theta-\pi $ to $ -\theta+\pi$), and finally a straight line from $\rho$ to  $\rho + e^{i(-\theta+\pi)}\infty$  (on the other bank of the cut). In the same cut plane, for $R>|\rho|$ we denote by $\Gamma'_R$  the circle minus one point defined by $|z| = R$ and $-\theta-\pi < \arg (z-\rho)  < -\theta+\pi$, positively oriented.

We can now state the relation with the inverse Laplace transform. Given $\theta\in\R$, we see any $f\in\Nilsrho$ as  a function holomorphic on 
\begin{equation} \label{petitdisque}
\{z\in\C, \, \, |z-\rho| < \eps, \,   \, -\theta-\pi < \arg (z-\rho)  < -\theta+\pi\}
\end{equation}
for some $\eps>0$ by letting $\log   (z-\rho) = \log|  z-\rho | + i \arg  (z-\rho) $.

\begin{prop} \label{proplienlaplaceinverserho}
Let $\rho\in\C$, $\theta\in\R$, and $f\in\Nilsrho$. Assume that $f$, seen as a function on \eqref{petitdisque}, can be analytically continued (for some $\eps > 0$) to both
$$\{z\in\C, \, \, |z| > |\rho| \mbox{ and } -\theta-\pi < \arg(z-\rho) < -\theta-\pi+\eps\}$$
and
$$\{z\in\C, \, \, |z| > |\rho| \mbox{ and } -\theta+\pi -\eps < \arg(z-\rho) < -\theta+\pi \},$$
with sub-exponential growth on these sectors. 
Then $\Cinf \Lrho f $ is the asymptotic expansion of $\frac1{2i\pi} \int_{\Gamma_\rho} f(z) e^{xz} \dd z$ as $|x|\rightarrow \infty$ in a large sector bisected by $\theta$.
\end{prop}

 By sub-exponential growth on these sectors $U$, we mean that for any $\epsz>0$ there exists $c_\epsz>0$ such that, for any $z\in U$ with $|z| \geq |\rho|+1$, we have  
$|f(z)| \leq c_\epsz \exp(\epsz |z|)$. 

\bigskip

On the other hand, given   $\theta\in\R$ any $f\in\Nilsinf$ yields a function holomorphic on 
\begin{equation} \label{domaineholonilsinf}
\{z\in\C, \, \, |z | > R_0 ,  \, -\theta-\pi < \arg(z)  < -\theta+\pi \}
\end{equation}
provided $R_0$ is sufficiently large, by letting $\log  z = \log|z| + i \arg(z)$. 

\begin{prop} \label{proplienlaplaceinverseinf}
Let $\theta\in\R$ and $f\in\Nilsinf$. Then in the expansion \eqref{eqdeferho} of $\Czero \Linf f \in \Ezero$, the function $\sum_{n=0}^\infty  c_{\alpha, j, n} x^n$ is entire for any $\alpha \in S$ and any $j\in T$. Moreover for $x\neq 0$ with $ \theta-\frac{\pi}{2} < \arg(x)  <  \theta+\frac{\pi}{2}$ we have
$$(\Czero \Linf f )(x)=   \lim_{R\rightarrow+\infty} \frac1{2i\pi}\int_{\Gamma'_R} f(z) e^{xz} \dd z$$
where $f$ is seen in the right hand side as a function holomorphic on \eqref{domaineholonilsinf}.
\end{prop}

Finally, we have the following result related to the usual Laplace transform. It is in disguise nothing but 
a form of the classical Watson's lemma.

\begin{prop} \label{prop7} Let  $\theta\in\R$, and $\eps>0$. Let $U$ denote the sector defined by $x\neq 0$ and $\theta-\eps < 
\arg(x) < \theta+\eps$. Let $g$ be a function holomorphic on $U$ such that $|g(x) | \leq A \exp(B |x|)$ for any $x\in U$, where $A,B>0$ are fixed. Assume also that $g$ has a generalized Taylor expansion at 0 in the Nilsson class $\Nilszero$, and that $g$ is locally integrable at 0.

Then the Laplace transform $\int_0^{e^{i\theta}\infty} g(x)e^{-zx} \dd x$ is defined and holomorphic for $z$ in a large sector bisected by $-\theta$ with $|z|$ sufficiently large, and its asymptotic expansion as $|z| \rightarrow \infty$ in this sector is $\Linf^{-1}  \Czero^{-1} g $ (where $g$ is seen in $\Nilszero$). 
\end{prop}
In this statement we recall that $\Linf : \Einf\rightarrow\Ezero$ is not bijective: its kernel is $\cz$ (see Proposition \ref{propker}). However it induces by restriction a bijective linear map $E'_{\infty} \rightarrow E'_0$, where $E'_{\infty}$ (resp.  $E'_0$) is the space of formal functions \eqref{eqdefeinf} (resp. \eqref{eqdeferho} with $\rho=0$) with $S \subset \{\alpha\in\C, \, \Re \alpha > -1\}$.

\bigskip

We omit the proof of Proposition \ref{prop7}, which is given in~\cite[p. 121, (4.4.17)]{BH} in another form: to 
make the connection, it is enough to observe that  $\Linf^{-1}\circ \Czero^{-1}$ maps $x^{-\alpha-1}(\log(1/x))^i$ to
\begin{equation}\label{eqprollap1}
\Big(\frac{\dd}{\dd y}\Big)^i\Big[ \Gamma(-y)z^y\Big]_{|y=\alpha}
 = \sum_{k=0}^i \binom{i}{k} (-1)^{i-k}\Gamma^{(i-k)}(-\alpha)z^\alpha (\log(z))^k  
 \end{equation}
provided $ \alpha\not\in\N$.

Let us now prove Propositions \ref{proplienlaplaceinverserho} and \ref{proplienlaplaceinverseinf}.
 Proposition \ref{proplienlaplaceinverserho} is a generalization of the claim made in the proof of Theorem 6 of \cite{ateo} and can be proved along the same lines; the only new ingredient is the following computation, valid   for $\alpha\in\C$ and $j \in\N$:
\begin{align*}
\Cinf &\Lrho \Big((z-\rho)^\alpha (\log(z-\rho))^j\Big) \\
&=  e^{\rho x} \sum_{k=0}^j \binom{j}{k} \Big( \frac{\Gamma(1-\{y\})}{\Gamma(-y)}\Big)^{(j-k)}(\alpha) \Cinf \Big( x^{-\alpha-1} (\log(1/x))^k\Big)\\
&=e^{\rho x}  x^{-\alpha-1}  \sum_{k=0}^j \binom{j}{k} \Big( \frac{\Gamma(1-\{y\})}{\Gamma(-y)}\Big)^{(j-k)}(\alpha) \sum_{i=0}^k (-1)^{k-i} 
\binom{k}{i} \Gh^{(k-i)}(1-\{\alpha\})(\log(1/x))^i \\
&=e^{\rho x}  x^{-\alpha-1}  \sum_{i=0}^j \binom{j}{i} (\log(1/x))^i \sum_{\ell = 0 }^{j-i} \binom{j-i}{\ell}   \Big( \frac{\Gamma(1-\{y\})}{\Gamma(-y)}\Big)^{(j-i-\ell)}(\alpha) ( \Gh(1-\{y\}))^{(\ell)} (\alpha) \\
&=e^{\rho x}  x^{-\alpha-1}  \sum_{i=0}^j \binom{j}{i} (\log(1/x))^i (\Gh(-y))^{(j-i)}(\alpha)\\
&=e^{\rho x} \Big(\frac{\dd}{\dd y}\Big)^{j} \Big[ \frac{x^{-y-1}}{\Gamma(-y)}\Big] _{|y = \alpha}.
\end{align*}

\bigskip

To conclude, let us prove Proposition  \ref{proplienlaplaceinverseinf}. Let  $\theta\in\R$, $f\in\Nilsinf$, and $R_0$  be sufficiently large.  Then $f(z)$ is given on \eqref{domaineholonilsinf}  by  a  convergent expansion 
$$
f(z) = \sum_{\al\in S} \sum_{j\in T} \sum_{n=0}^{\infty} c_{\al,j,n}z^{-n-\al-1}(\log(1/z))^j
$$
where $S\subset\C$ and $T\subset\N$ are finite subsets, and $c_{\al,j,n} \in\C$.  
We then have
\begin{align*}
\frac1{2i\pi} \int_{\Gamma'_R} f(z)e^{zx}\dd z &= \sum_{\al\in S} \sum_{j\in T} \sum_{n=0}^{\infty} c_{\al,j,n} 
\frac1{2i\pi} \int_{\Gamma'_R} z^{-n-\al-1}(\log(1/z))^j e^{zx} \dd z\\
&=\sum_{\al\in S} \sum_{j\in T} \sum_{n=0}^{\infty} c_{\al,j,n} \frac{\partial^j}{\partial \alpha^j} \left[
\frac1{2i\pi} \int_{\Gamma'_R} z^{-n-\al-1} e^{zx} \dd z\right];
\end{align*}
recall that $ -\theta-\pi < \arg(z)  < -\theta+\pi$, $\log(1/z) = -\log(z)$, and the cut corresponds to $\arg(z) = -\theta \pm \pi$. Assume that $ \theta-\frac{\pi}{2} < \arg(x)  <  \theta+\frac{\pi}{2}$; then $\arg(zx)$ belongs to either $(\frac{\pi}{2}, \frac{3\pi}{2})$ or $(-\frac{3\pi}{2},  -\frac{\pi}{2}) $ on each bank of the cut.  If $\Re(-n-\alpha)>0$, we can flatten $\Gamma'_R$ on the cut and let $R\rightarrow +\infty$; a simple 
computation then shows that 
\begin{equation}\label{eq:intgammaR}
\lim_{R\rightarrow +\infty}\frac1{2i\pi} \int_{\Gamma'_R} z^{-n-\al-1} e^{zx} \dd z = \frac{x^{n+\alpha}}{\Gamma(n+\alpha+1)}.
\end{equation}
If $\Re(-n-\alpha)\le 0$, we integrate enough times by parts (all the integrated parts vanish) to come back to the above situation; in the end,~\eqref{eq:intgammaR} holds again.

We deduce that when $ \theta-\frac{\pi}{2} < \arg(x)  <  \theta+\frac{\pi}{2}$, 
\begin{equation} \label{eqLapinfsomme}
\lim_{R\rightarrow +\infty}\frac1{2i\pi} \int_{\Gamma'_R}  f(z)e^{zx}\dd z 
 = \sum_{\al\in S} \sum_{j\in T} 
\sum_{n=0}^{\infty} c_{\al,j,n} \frac{\partial^j}{\partial \alpha^j} \left[\frac{x^{n+\alpha}}{\Gamma(n+\alpha+1)}
\right].
\end{equation}
Letting $x$ tend to 0 and noticing that   the above computation of $\Cinf \Lrho \Big((z-\rho)^\alpha (\log(z-\rho))^j\Big)$ yields with $\rho=0$:
$$\Czero \Linf  \Big( z^{-n-\alpha-1}  (\log(1/z))^j\Big) = (-1)^j \Big(\frac{\dd}{\dd y}\Big)^j \Big[ \frac{x^{-y-1}}{\Gamma(-y)}\Big]_{|y=-n-\al-1},$$
this concludes the proof of Proposition  \ref{proplienlaplaceinverseinf}.

\bigskip

\begin{Remark} \label{rem2}
As Proposition \ref{prop7} shows, the map $\Linf^{-1}\circ \Czero^{-1}$ extends the usual Laplace transform, because it is formal, so that no convergence assumption is needed. It maps $x^{-\alpha-1}(\log(1/x))^i$ to 
$$\frac1{i+1} \sum_{k=0}^{i+1} \binom{i+1}{k}  \Big(\frac{\dd}{\dd y}\Big)^{i+1-k} \Big( (y-\alpha)\Gamma(-y)\Big)_{|y=\alpha}
 z^\alpha (\log(z))^k  \, \, \bmod  \C[[z]]$$
for any $\alpha\in\C$ and any $i\in\N$;  if $ \alpha\not\in\N$ this formula can also be written as \eqref{eqprollap1}.
  These are exactly the implicit formulas (5.3.10) obtained by Andr\'e  \cite{YA1} using operational calculus. Therefore one may use 
them instead of operational calculus in Andr\'e's proof. However our approach in \S \ref{subsecnvpreuve} is different: we use only $\Lseul$ so that no transcendental coefficient appears.
  \end{Remark}

\section{Microsolutions, analytic continuation, and  Laplace transform} \label{sec3}

In this section we prove Theorem \ref{thintroun} announced in the introduction (see \S \ref{subsec32}), after restating it in \S \ref{secmicro} in terms of microsolutions and analytic continuation. Assuming $\calD$ to be Fuchsian, we relate this theorem in \S \ref{subseclaplace} to Laplace transform, and then in \S \ref{secdiagLG} to a duality between the solutions of $\calFbar\calD$ at 0 and at $\infty$; we use it to obtain a new construction of Andr\'e's extension of Laplace transform. We conclude this section with a commutative diagram which summarizes the results of \S \S \ref{sec2} and \ref{sec3}, including the factorization of the inverse Laplace transform obtained in \S \ref{subsec25}.

\subsection{Setting and notation} \label{secmicro}

Let $\calD \in \C[z, \frac{\dd}{\dd z}]$ be a differential operator of which $\infty$ is a regular singularity (or even not a singularity at all). We  denote by $\Sigma$ the set of all finite singularities of $\calD$; some of them might be apparent singularities but play a crucial role in our setting: see Remark \ref{remapp} in \S \ref{subsec44}.  We also fix a real number $\theta$ such that $\theta\not\equiv -  \arg(\rho-\rho') \bmod \pi$ whenever $\rho,  \rho' \in \Sigma$ are distinct. We shall work   in   the simply connected  cut plane $\cut $ obtained from $\C$ by removing the union of all closed half-lines of direction $- \theta + \pi $ starting at elements of $\Sigma$. For $z$ in this cut plane and $\rho\in\Sigma$, we agree that $-\theta-\pi < \arg(z-\rho) < -\theta+\pi$. 

\bigskip

Let us denote by $\delta$ the degree of $\calD$, by $\mu$ its order, and write
$$\calD  = \sum_{j=0}^\mu P_{\mu-j} (z) (\frac{\dd}{\dd z})^j$$
where $P_i(z) \in\cz$ has degree at most $\delta-i$ (with equality for $i=0$), since $\infty$ is a regular singularity. For any $j\in\zeromu$ let $a_{\mu-j} $ denote the coefficient of degree $\delta-\mu+j$ of $P_{\mu-j}(z)$, so that the indicial equation at infinity is $R(-z) = 0$, where
\begin{equation} \label{eqdefR}
R(z) = \sum_{j=0}^\mu a_{\mu-j} z(z-1)\ldots (z-j+1) \in \cz\setminus\{0\}.
\end{equation}
Let $e_1,\ldots,e_p$ denote the non-negative  integer roots of $R$ (with $p=0$ if there is no such root), so that $-e_1$, \ldots, $-e_p$ are the non-positive integer exponents of $\calD$ at infinity.  We let
$$M = \max(0, \mu-\delta, e_1,\ldots, e_p) \hspace{0.8cm} \mbox{ and }  \hspace{0.8cm}  N = M + \delta-\mu \geq 0.$$
We shall prove that in  Theorem \ref{thintroun} the function $f$ can be chosen such that $\deg\calD f < N$.

\bigskip

For any $\rho\in\Sigma$, recall that $\Orho$ and $\Orhoch$ are defined in the introduction, and let $\Srho$ denote the kernel of $\calD$ seen as a linear map $\Orhoch/\Orho \rightarrow \Orhoch/\Orho $. In other words, 
$$
\Srho  =
\{f\in\Orhoch, \, \, \calD f \in \Orho\}/\Orho.
$$
This is the space of microsolutions of $\calD$ at $\rho$. Kashiwara's theorem \cite{Kashiwara}  (see also \S 1.2 of \cite{Pham})
 asserts that $\dim \Srho =  \mrho$,    the multiplicity of $\rho$ as a singularity of $\calD$. We also let 
$$\Sinf =  \ker \Big(   \frac{\hol }{\cz}   \stackrel{\calD}{\rightarrow}   \frac{\hol }{\cz} \Big) =   \frac{\calD^{-1}(\cz)}{\cz}$$
which is {\em not}   the space of microsolutions at infinity; here $\hol$ is  the space of holomorphic functions on $\cut$ and   $\calD^{-1}(\cz) = \{ f\in\hol,\, \calD f \in\cz\}$. 

\bigskip

Now any $f\in \hol$ can be restricted to a small  cut disk \eqref{eqdisquecoupe} around any given $\rho\in\Sigma$. This provides a map $\hol\rightarrow\Orhoch$, which induces a linear map $\Sinf\rightarrow\Srho$ and then a diagonal map $\pro : \Sinf \rightarrow \oplus_{\rho \in \Sigma} \Srho$ by mapping $f$ to $(f,f,\ldots,f)$, where the coordinate $f$ corresponding to $\rho\in\Sigma$ is seen locally around $\rho$. This map $\pro$  can be thought of as analytic continuation in $\cut$ from a neighborhood of  infinity towards a neighborhood of all singularities $\rho\in\Sigma$.

\bigskip

With these notations, Theorem \ref{thintroun} is equivalent to the following result.

\begin{theo}\label{thmicrosol} The map $\pro : \Sinf \rightarrow \oplus_{\rho \in \Sigma} \Srho$ is bijective.
\end{theo}

We shall also prove that $\Sinf =  \frac{\calD^{-1}(\czN)}{\cz\cap \calD^{-1}(\czN)} $; this shows that in Theorem \ref{thintroun} there exists $f$ such that $\deg\calD f <N$. This implies 
\begin{equation} \label{eqcaractsrho}
\Srho = \frac{\calD^{-1}(\czN)}{\Orho\cap\calD^{-1}(\czN)}
\end{equation}
where   $\calD^{-1}(\czN) = \{ f\in\hol,\, \calD f \in\czN\}$. This equality will be used in \S \ref{subsecnvpreuve} because when $\calD$ is a $G$-operator, $\calD^{-1}(\czN)$ is the space of solutions in $\hol$ of $(\frac{\dd}{\dd z})^N\circ\calD$ which is also a $G$-operator; therefore the Andr\'e-Chudnovski-Katz Theorem can be applied.

\medskip

The special case $N=0$ is also worth mentioning,  because $\calD^{-1}(\cz_{<0}) = \ker\calD$.

 \begin{coro} \label{cormicro}
 Let  $\calD \in \C[z, \frac{\dd}{\dd z}]$ be a differential operator of which $\infty$ is a regular singularity; assume that 
  $\mu\geq\delta$ and no integer exponent of $\calD$ at $\infty$ is less than $ \delta - \mu$. Then
 in Theorem~\ref{thintroun} there exists $f$ such that $\calD f=0$, and  for any finite singularity $\rho$:
 \begin{itemize}
 \item[$(i)$] 
Any microsolution in $\Srho$ can be represented by a solution of $\calD$.
 \item[$(ii)$]  We have $\calD^{-1}(\Orho) = \Orho + \ker \calD$, and $\calD$ induces a surjective map $\Orho\rightarrow\Orho$. 
\end{itemize}
\end{coro}

Indeed  given $g\in\Orhoch$ such that $\calD g\in\Orho$, we find $f\in\ker\calD$  which has the same class as $g$ in $\Srho$: this proves $(i)$. To prove $(ii)$, given $h \in \Orho$ we consider $g\in\Orhoch$ such that $\calD g=h$, and we have  $\calD (g-f)=h$ with $g-f\in\Orho$.

\subsection{Proof  of  Theorem \ref{thmicrosol}} \label{subsec32}

The proof of  Theorem \ref{thmicrosol} splits into two parts.

\bigskip

The first step is to prove  that, as a $\C$-vector space,
\begin{equation} \label{eqlemdim}
\Sinf =  \frac{\calD^{-1}(\czN)}{\cz\cap \calD^{-1}(\czN)} \mbox{ has dimension } \delta = \deg \calD;
\end{equation}
here   $\calD^{-1}(\czN) = \{ f\in\hol,\, \calD f \in\czN\}$. Taking this for granted, let us conclude the  proof of  Theorem \ref{thmicrosol}.

 Let $f\in\calD^{-1}(\cz)$ be such that $f+ \cz \in \ker\pro$. Then $(\frac{\dd}{\dd z})^{\deg \calD f + 1} \calD f = 0$,  and the finite singularities of $(\frac{\dd}{\dd z})^{\deg \calD f + 1} \circ \calD$ are the same as those of $\calD$, so that $f$ is holomorphic at any point outside $\Sigma$. Now $f+ \cz \in \ker\pro$ so that $f$ is also holomorphic at any $\rho\in\Sigma$. Therefore $f$   is entire, and using a generalization of a theorem of Liouville (see for instance   \cite{Ahlfors}),   we deduce that $f\in\cz$. Therefore 
  $\pro$ is injective.

Now combining \eqref{eqlemdim} with Kashiwara's theorem, we obtain
$$\dim\Sinf = \delta = \sum_{\rho\in\Sigma} \mrho = \sum_{\rho\in\Sigma} \dim \Srho
$$
since $\deg\calD$ is the degree of the coefficient $P_0(z)$ of $(\frac{\dd}{\dd z})^\mu$ in $\calD$  because  $\infty$ is a regular singularity. Therefore $\pro$ is bijective; this concludes  the proof of Theorem \ref{thmicrosol}.

\bigskip

\begin{Remark} \label{rem3}
Using \eqref{eqlemdim} we obtain $f$ in Theorem \ref{thintroun} such that $\deg\calD f < N$. 
\end{Remark}

Let us prove Eq. \eqref{eqlemdim} now. The proof is a variant of that of Th\'eor\`eme 1.4 of \cite{MalgrangePointsSing}. The main step is to compute the index of $\calD : \cz\rightarrow\cz$ (see Eq. \eqref{eqdimdinf} below), as a special case of   the proof of Th\'eor\`eme 2 of \cite{BertrandKummer}.

For any $k \geq \max(0, \mu-\delta)$ we have
$$\calD z^k = R(k) z^{k+\delta-\mu} + Q_k(z),$$
where $Q_k(z)\in\cz $ has degree less than $k+\delta-\mu$ and $R(z)$ is the polynomial \eqref{eqdefR}. Therefore $\calD$ induces, for any $n\geq  \max(0, \mu-\delta)$, a linear map
$$\calD_n : \cz_{<n} \rightarrow \cz_{<n+\delta-\mu}$$
where $\cz_{<0} = \{0\}$. Since $\calD z^k$ has degree exactly $k+\delta-\mu$ for any $k\geq M$ (by definition of $M$), we have
$$\calD(\cz) = \calD(\czM) \oplus \Span\{z^j, \, j \geq N\}.$$
Therefore the identity map induces a bijective linear map 
\begin{equation} \label{eqdimdn}
\frac{\czN}{\calD(\czM)}  \stackrel{\sim}{\longrightarrow} \frac{\cz }{\calD(\cz )}
\end{equation}
so that
\begin{equation} \label{eqdimdinf}
\dim\Big( \frac{\cz }{\calD(\cz )}\Big) = N - \rk \calD_M = \delta- \mu + \dim\ker\calD_M .
\end{equation}
In other words, the index of $\calD : \cz\rightarrow\cz$ is $\mu-\delta$ (see \cite{MalgrangePointsSing}).

\bigskip

To conclude the proof of \eqref{eqlemdim}, one may use exact sequences as in the proof of Th\'eor\`eme 1.4 of    \cite{MalgrangePointsSing}; for the convenience of the reader we provide a more down-to-earth proof.
Let $t = M$ or $t= \infty$, with $\cz_{< \infty} = \cz$. We consider  the linear map
$$\widetilde \calD : \frac{\calD^{-1} (\cz_{< t + \delta-\mu})}{\cz \cap \calD^{-1} (\cz_{< t + \delta-\mu})} \rightarrow
 \frac{\cz_{< t + \delta-\mu} }{\calD(\cz _{< t })}$$
induced by $\calD$; here $\calD^{-1} (\cz_{< t + \delta-\mu})$ is the set of all $f\in\hol$ such that $\calD f \in \cz_{< t + \delta-\mu}$. This linear map
is surjective since $\calD : \hol\rightarrow\hol$ is surjective; its kernel is
$$\ker \widetilde \calD  = \frac{\calD^{-1} (\calD( \cz_{< t  }))}{\cz \cap \calD^{-1} (\calD( \cz_{< t  }))} = \frac{\cz_{<t}+\ker\calD}{\cz_{<t}} = 
\frac{\ker\calD}{\cz_{<t} \cap \ker\calD}$$
where $\ker\calD = \{f\in\hol, \, \calD f =0\}$ has dimension $\mu$. Therefore we obtain, using Eqns. \eqref{eqdimdn} and  \eqref{eqdimdinf} and the assumption $t\in\{M, \infty\}$:
\begin{align*}
\dim \frac{ \calD^{-1} (\cz_{< t + \delta-\mu})}{\cz \cap \calD^{-1} (\cz_{< t + \delta-\mu})} 
&=  \dim\ker\widetilde \calD + \rk \widetilde \calD \\
&=  \mu - \dim(\cz_{<t} \cap \ker\calD) + \dim\Big(\frac{  \cz_{< t + \delta-\mu}}{\calD(\cz_{<t})} \Big) 
\\
&=  \mu - \dim\ker\calD_M + \dim\Big(  \frac{\cz }{\calD(\cz )} \Big) = \delta.
\end{align*}
With $t=\infty$ this means  $\dim\Sinf = \delta$; with $t=M$ we deduce that   the natural injective map
$$  \frac{ \calD^{-1} (\czN)}{\cz \cap \calD^{-1} (\czN)}  \rightarrow   \frac{ \calD^{-1} (\cz )}{\cz  }  $$
induced by the identity map is an isomorphism. This concludes the proof of \eqref{eqlemdim}.

\bigskip

\begin{Remark}\label{rem4}
 In this proof $\hol$ may be replaced with $\Nilsinf$ since $\infty$ is a regular singularity. We obtain in this way:
\begin{equation}\label{eqegalSinf}
\Sinf =  \ker \Big(   \frac{\Nilsinf}{\cz}   \stackrel{\calD}{\longrightarrow}   \frac{\Nilsinf }{\cz} \Big) =   \frac{\{f\in\Nilsinf,\, \, \calD f \in\cz\}}{\cz}
= \frac{\{f\in\Nilsinf,\, \, \calD f \in\czN\}}{\{f\in\cz ,\, \, \calD f \in\czN\}}.
 \end{equation}
\end{Remark}

\subsection{Laplace transforms} \label{subseclaplace}

We keep the notation and assumptions of \S \ref{secmicro}, and we assume (from now on) that $\calD$ is Fuchsian.
This subsection is essentially an adaptation to our context of \S 1.5 of \cite{Pham}.

We denote by $\calFbar : \C[z,\frac{\dd}{\dd z}] \rightarrow \C[x,\frac{\dd}{\dd x}]$ the Fourier transform of differential operators, i.e. the morphism of $\C$-algebras defined by $\calFbar(z) =   \frac{\dd}{\dd x}$ and $\calFbar(\frac{\dd}{\dd z}) = -x$.

We denote by $\Sbarth$ the space of solutions of $ \calFbar \calD $  holomorphic on the simply connected cut plane $\cutbar $ defined by $x\neq 0$ and $\theta-\pi < \arg(x)  <  \theta+\pi$. It has dimension $\delta$ (i.e., equal to the order of $ \calFbar \calD $ which is the degree of $   \calD $) since $ \calFbar \calD $ has no singularity except 0 and $\infty$ (the proof \cite{YA1} of this property in the special case where $\calD$ is a $G$-operator applies to any Fuchsian operator). This space can be written (see \cite{Bresil} or \cite{Malgrangelivre}) as $\Sbarth = \oplus_{\rho\in\Sigma} \Sbarthrho$ where $\Sbarthrho$ is the set of $g\in\Sbarth$ such that $g(x) e^{-\rho x}$ has sub-exponential growth as $|x|\rightarrow+\infty$ in a large sector bisected by $\theta$ (i.e., such that for any $\eps>0$ there exists $c_\eps>0$ with $|g(x)e^{-\rho x}| \leq c_\eps \exp(\eps |x|)$ for any $x$ in   a large sector bisected by $\theta$ with $|x|\geq 1$).

\bigskip

A local Laplace transform  has been studied by Malgrange and Ecalle. Given $\rho\in\Sigma$, they prove that letting
\begin{equation} \label{eqdeflaprho}
(\Vrho g)(z) =  \int_{x_0} ^{e^{i\theta}\infty}  g(x) e^{-zx} \dd x \, \bmod \Orho
 \end{equation}
for $g\in  \Sbarthrho$ and $z$ sufficiently close to $\rho$ with $-\theta-\frac{\pi}{2} < \arg(z-\rho) < -\theta+\frac{\pi}{2}$ (where $x_0 \in \cutbar$ is chosen arbitrarily) provides a linear isomorphism
$$\Vrho : \Sbarthrho   \stackrel{\sim}{\longrightarrow} \Srho  
=
\{f\in\Orhoch, \, \, \calD f \in \Orho\}/\Orho 
$$
(see Theorem 2.2 of \cite{Bresil}). For $x_0,x'_0\in\cutbar$ the map $z\mapsto   \int_{x_0} ^{x'_0}  g(x) e^{-zx} \dd x $ is holomorphic at $\rho$, so that $ \Vrho g  \in \Srho$ is independent from the choice of $x_0$.
The inverse image of the class of  $f\in\Orhoch$ such that $\calD f \in \Orho$ is given by
\begin{equation} \label{eqdefgtube}
(\Vrho^{-1}(f\bmod \Orho)) (x) = \frac1{2i\pi}\int_{\Gamma_\rho} f(z) e^{xz} \dd z
\end{equation}
for $\theta-\frac{\pi}{2} < \arg (x)  < \theta+\frac{\pi}{2}$, 
where $\Gamma_\rho$ is  the following path in the cut plane $\cut$ (as  in \S \ref{subsec25}): a straight line from $\rho + e^{i(-\theta-\pi)}\infty$ to $\rho$ (on one bank of the cut), a circle of radius essentially zero around $\rho$ (with $\arg (z-\rho)  $ increasing from $-\theta-\pi $ to $ -\theta+\pi$), and finally a straight line from $\rho$ to  $\rho + e^{i(-\theta+\pi)}\infty$  (on the other bank of the cut). Here we use the fact that $f$ can be analytically continued to both sectors of Proposition \ref{proplienlaplaceinverserho}; actually $f$ can be chosen to be holomorphic on $\cut$ using Theorem \ref{thintroun}, and it has moderated growth at infinity. Note also  that $x$ has been changed into $-x$ with respect to \cite{Bresil} since we consider here $ \calFbar$ instead of  $\calF = \calFbar^{-1}$.

\bigskip

These ``local'' Laplace transforms at various $\rho\in\Sigma$ can be glued together to obtain a linear isomorphism
$$\Vloc  : 
 \Sbarth =   \oplus_{\rho\in\Sigma} \Sbarthrho  \stackrel{\sim}{\longrightarrow}\oplus_{\rho\in\Sigma} \Srho  $$
by letting $\Vloc ( \sum_\rho g_\rho) = \sum_\rho \Vrho (g_\rho)$.
Now Theorem \ref{thmicrosol} provides a linear isomorphism $\pro : \Sinf\rightarrow\oplus_\rho \Srho$. Therefore we can define a bijective linear map
\begin{equation} \label{eqdefLapprolong}
\Vinf = \pro^{-1} \circ \Vloc :  \Sbarth =   \oplus_{\rho\in\Sigma} \Sbarthrho  \stackrel{\sim}{\longrightarrow} \Sinf.
\end{equation}
We shall prove at the end of \S \ref{secdiagLG} below that for any $g\in  \Sbarth$ locally integrable around 0, we have
\begin{equation} \label{eqdefVinf}
(\Vinf g)(z) = \int_0^{e^{i\theta}\infty} g(x) e^{-zx} \dd x
\end{equation}
for any $z\in\C$ such that $|z|$ is large enough and $-\theta-\frac{\pi}{2} < \arg  z < -\theta+\frac{\pi}{2}$. Therefore $\Vinf$ extends the usual Laplace transform (taking $\theta =0$, say) to any $g\in  \Sbarth$: we recover Andr\'e's extension \cite{YA1} based on operational calculus. See 
also Remark \ref{rem2} at the end of \S \ref{subsec25}. 

\bigskip

To compute the inverse map $\Vinf^{-1}$ we   have to define  another path of integration (see    \cite[pp. 183--192]{DiP} or \cite{ateo}). Let $R$ be such that $R > |\rho| $ for any $\rho\in\Sigma$. 
Given $\rho \in \nvsing$, let  $\deltarho = \rho - e^{-i\theta}\R_+$ denote  the half-line of 
angle $-\theta+\pi \bmod 2\pi$ starting at $\rho$. It  intersects
  the circle $\cercle $  centered at 0 of radius $R$  at one point 
$\pointrho = \rho - \Arho e^{ -i\theta}$, with $\Arho >0$, which corresponds to two points 
at the border of the cut plane $\cut$, namely $\rho + \Arho e^{i(-\theta\pm\pi)}$ with values $-\theta\pm\pi$ 
of the argument. We denote by  $\Gamma'_R$  the path   going in the positive direction  from $z_{\rho_1}$ to  $z_{\rho_p}$ 
along the circle ${\mathcal C}(0,R)$. Here $\rho_1,\rho_p \in\Sigma$ are such that the minimal (resp. maximal) value of $\Im (\rho e^{i\theta})$ for $\rho\in\Sigma$ is taken at $\rho = \rho_1$ (resp. $\rho = \rho_p$). In this way, when $R$ is very large, $\Gamma'_R$ is nearly the whole circle and this is essentially the same notation $\Gamma'_R$ as in  \S \ref{subsec25}.

\begin{figure} 
\centering
\includegraphics[scale=0.5, width=0.5\textwidth ,trim = 0 320 0 80, clip=true]{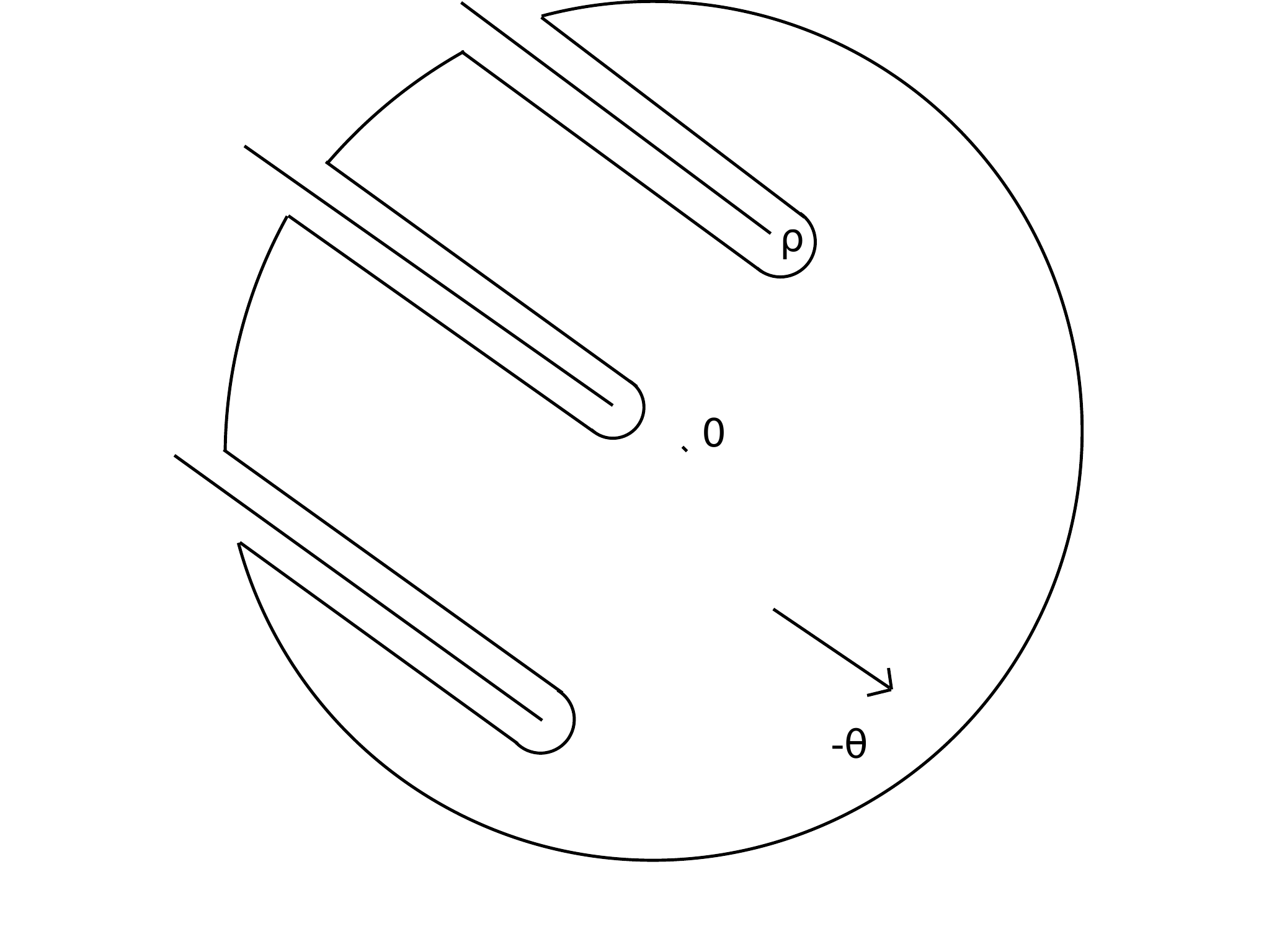}
\caption{The contour $\Gamma_R$} \label{fig1}
\end{figure}

Now let $f\in\hol$ be such that $\calD f \in \cz$. Then we have
\begin{equation} \label{eqlapinfmu}
 (\Vinf^{-1}(f \bmod \cz) )(x)=  \lim_{R\rightarrow+\infty} \frac1{2i\pi}\int_{\Gamma'_R} f(z) e^{xz} \dd z 
 \end{equation}
for any $x$ such that $ \theta-\frac{\pi}{2} < \arg  x <  \theta+\frac{\pi}{2}$. Indeed $\Vinf^{-1} f = \sum_{\rho\in\Sigma} \Vrho^{-1} f_\rho$ where $f_\rho$ is the image of $f$ in $\Srho$, so that Eq. \eqref{eqlapinfmu} means
$$  \sum_{\rho\in\Sigma}  \frac1{2i\pi}\int_{\Gamma_\rho} f(z) e^{xz} \dd z = \lim_{R\rightarrow+\infty} \frac1{2i\pi}\int_{\Gamma'_R} f(z) e^{xz} \dd z .$$
This equality follows from  the residue theorem applied on the contour of Fig. \ref{fig1} by letting $R\rightarrow\infty$. 

\begin{Remark} \label{rem5}
These Laplace transforms provide  a way to compute $\pro^{-1}$, that is, given microsolutions $f_\rho\in\Srho$ for any $\rho\in\Sigma$, to find an analytic  solution $f$ of $\calD$ on the cut plane $\cut$ such that, for any $\rho$, the function $f-f_\rho$ is holomorphic at $\rho$. Indeed, one may take
$$f = \sum_{\rho\in\Sigma} \Vinf  \Invlaprho f_\rho $$
and use Eq. \eqref{eqdefVinf} if it applies, or the operators $\Linf^{-1}$ and $\Czero^{-1}$ (see below), to compute $f$.
\end{Remark}

\subsection{Relation with the formal operators $\Lseul$ and $\Cseul$}\label{secdiagLG}

In this section we keep the assumptions and notation of \S \ref{subseclaplace}, and study the relation with the formal operators of \S \ref{sec2}. We shall summarize most constructions of \S \S \ref{sec2} and \ref{sec3} in a commutative diagram.

 For any $\rho\in\Sigma$ we denote by $ \Sbarinfrho$ the space of formal solutions at infinity of $\calFbar \calD$ in $e^{\rho x} \Einf$; it has dimension $\mrho$. See \cite{YA1}, where  the proof extends to any Fuchsian operator.   In general the elements of $\Einf$ that appear here involve  divergent series.  

Taking asymptotic expansion as $|x|\rightarrow+\infty$ in a large sector bisected by $\theta $ (see \cite{Ramis} or  \cite{ateo} for the precise definition) yields a linear isomorphism 
$$\asy : \Sbarthrho  \stackrel{\sim}{\longrightarrow}  \Sbarinfrho,$$
the inverse of which is Borel-Laplace summation $ \asymu$ (i.e., Ramis' 1-summation) since $\theta$ is not anti-Stokes (recall that $\theta\not\equiv -  \arg(\rho-\rho') \bmod \pi$ whenever $\rho,  \rho' \in \Sigma$ are distinct). This extends by linearity to a linear isomorphism still denoted by the same symbol:
$$\asy :\Sbarth = \oplus_{\rho\in\Sigma}  \Sbarthrho  \stackrel{\sim}{\longrightarrow}   \oplus_{\rho\in\Sigma}  \Sbarinfrho.$$

\bigskip

On the other hand, any $f\in\Sbarth$ is a holomorphic solution of $\calFbar \calD$ on the cut plane $\cutbar$ defined by $x\neq 0$ and  $\theta-\pi < \arg(x)  <  \theta+\pi$. Since 0 is a regular singularity of $\calFbar\calD$ (see \cite{YA1}), $f$ can be seen locally around the origin as an element of $\Ezero$. This provides a linear bijective map $\tay$ (namely, a generalized Taylor expansion)
$$\tay : \Sbarth   \stackrel{\sim}{\longrightarrow}   \Sbarzero$$
where $ \Sbarzero$ is the $\delta$-dimensional space of formal solutions of  $\calFbar\calD$ in $\Ezero$ (all of which involve series with a positive radius of convergence, so that $ \Sbarzero\subset\Nilszero$).

\bigskip

 Using also Propositions \ref{proplienlaplaceinverserho} and \ref{proplienlaplaceinverseinf} (proved in \S \ref{subsec25}), we obtain 
  that the following diagram is commutative:

$$
\bfig
\node a(0,300)[ \Sinf ] 
\node b(0,900)[  \oplus_\rho \Srho ] 
\node c(2000,600)[\Sbarth = \oplus_\rho \Sbarthrho] 
\node d(2000,1200)[\oplus_\rho \Sbarinfrho] 
\node e(2000,0)[\Sbarzero] 
\node f(1400,1200)[\oplus_\rho \Sbarinfrho] 
\node g(1400,0)[\Sbarzero] 
\arrow[a`b; \pro]
\arrow[a`c; \Invlapinf]
\arrow[b`c;   \Invlapth]
\arrow[c`d; \asy]
\arrow[c`e;\tay ]
\arrow[b`f;\oplus_\rho \Lrho]
\arrow[f`d;\Cinf]
\arrow[a`g;\Linf]
\arrow[g`e;  \Czero]
\efig 
$$

Here all maps are linear, bijective,  and depend on $\theta$. Only  the vector spaces $\Sbarinfrho$ and $\Sbarzero$ are independent from $\theta$ (for instance $\log(x)$ and $x^\alpha$, for $\alpha\in\C$, are formal quantities in $ \oplus_{\rho\in\Sigma}  \Sbarinfrho$ and $ \Sbarzero$).
To define $\oplus_\rho \Lrho$ and $ \Linf$, we identify $\Srho$ with a subspace of 
$
\Nilsrho/\C[[z-\rho]]$ using generalized Taylor expansion at $\rho$, and also $\Sinf $ with  a subspace 
of 
$
\Nilsinf/\cz$ using generalized Taylor expansion at infinity.

 \bigskip
 
 This diagram yields $\Vinf = \Linf^{-1} \circ \Czero ^{-1} \circ \tay$; using Proposition \ref{prop7} this proves Eq. \eqref{eqdefVinf} announced in \S \ref{subseclaplace}, namely that Eq. \eqref{eqdefLapprolong} provides a new construction of Andr\'e's extension of Laplace transform.

\section{Applications to special values} \label{secspecialvalues}

In this section we (mostly) assume that $\calD$ is a $G$-operator. To begin with, we 
recall Andr\'e's duality result on $E$- and \antie-functions, and give a new and constructive  proof of it built on our previous results. We also state and prove in \S \ref{subsecarith} our main result on the arithmetic nature of the coefficients of the matrices of all linear maps that appear in the commutative diagram of \S \ref{secdiagLG}. This result is then shown in \S \ref{subsec43}  to imply all arithmetic results stated in the introduction.  At last, we conclude this paper in \S \ref{subsec44} with an example related to Gompertz' constant.

\subsection{A new and constructive proof of Andr\'e's theorem} \label{subsecnvpreuve}

In this section we give a new and constructive proof of the following duality theorem (\cite{YA1}, Th\'eor\`eme 4.3).

\begin{theo}[Andr\'e] \label{thYA}
Let $\K$ be a number field, and $\calD \in \K[z, \frac{\dd}{\dd z}]$ be a $G$-operator of degree $\delta$; consider the differential equation $(\calFbar \calD)y=0$. Then:
\begin{itemize}
\item[$(i)$] There exists a basis of solutions around 0 of the form 
$$(F_1(x),\ldots,F_\delta(x)) \cdot x^{\Gamma_0}$$
where $F_1(x),\ldots,F_\delta(x)$ are $E$-functions with coefficients in $\K$,  and $\Gamma_0\in M_{\delta}(\Q)$ is an upper triangular matrix.
\item[$(ii)$] There exists a basis of formal solutions at infinity of the form 
$$(\ae_1(1/x),\ldots,\ae_\delta(1/x)) \cdot (1/x)^{\Gamma_\infty} \cdot e^{\Delta x}$$ 
where $\ae_1$, \ldots, $\ae_\delta$ are \antie-functions  with coefficients in the number field  $\K(\Sigma)$ generated by $\K$ and all finite singularities $\rho$ of $\calD$, $\Delta$ is the diagonal matrix with the finite singularities of $\calD$ as diagonal elements (repeated according to their multiplicities), and $\Gamma_\infty\in M_{\delta}(\Q)$ is an upper triangular matrix which commutes with $\Delta$.
\item[$(iii)$] The diagonal coefficients of $\Gamma_0$ (resp. $\Gamma_\infty$) are either in $\Z$, or congruent mod  $\Z$ to exponents of $\calD$ at infinity (resp. at finite singularities).
\end{itemize}
\end{theo}

\bigskip

In this section we deduce from the results of \S \S \ref{sec2}-\ref{sec3} a new proof of this result, leading to an effective construction of these bases. We 
also obtain a   generalization (see assertion $(i)$ below), upon noticing 
that everything works as soon as $\calD\in \K[z, \frac{\dd}{\dd z}] $ is Fuchsian with  exponents in $\K$ at all singularities, except that the $F_j$ and the $\ae_j$ are no  more $E$- or \antie-functions in general. We restrict to solutions of  $(\calFbar \calD)y=0$ at infinity, but the 
 same procedure can be carried out analogously to construct a basis of local solutions at 0 of $\calFbar \calD$ from a basis of $\Sinf$.

\bigskip

Let $\calD \in \C[z, \frac{\dd}{\dd z}]$ be a Fuchsian operator, and $N$ denote the integer defined in \S  \ref{secmicro}. 
We  denote by $\Sigma$ the set of all finite singularities of $\calD$, and we choose a real number $\theta$ such that $\theta\not\equiv -  \arg(\rho-\rho') \bmod \pi$ whenever $\rho,  \rho' \in \Sigma$ are distinct. As in \S \ref{sec3} we shall  work   in   the cut plane $\cut $ obtained from $\C$ by removing the union of all closed half-lines of direction $- \theta + \pi $ starting at elements of $\Sigma$. For $z\in\cut$ and $\rho\in\Sigma$, we agree that $-\theta-\pi < \arg(z-\rho) < -\theta+\pi$. 
For any $\rho\in\Sigma$ there   exist (non necessarily distinct) complex   numbers $t_1^\rho, \ldots, t_{J(\rho)}^\rho$, with $J(\rho)\geq 1$, and  functions $g_{j,k}^\rho(z-\rho)$ holomorphic at $\rho$, for $1\leq j \leq J(\rho) $ and $0\leq k \leq K(\rho,j)$,  such that the functions 
 \begin{equation}\label{eqdeffjk}
f_{j,k}^\rho(z-\rho) = (z-\rho)^{t_j^\rho} \sum_{k'=0}^k g_{j,k-k'}^\rho(z-\rho) \frac{(\log(z-\rho))^{k'}}{k'!}, 
\end{equation}
for $1\leq j \leq J(\rho) $ and $0\leq k \leq K(\rho,j)$, form a basis of solutions of $(\frac{\dd}{\dd z})^N\circ \calD$ in $\hol$ (where $\hol$ is the space of functions holomorphic on $\cut$). 
We  denote by $\Irho$ the set of pairs $(j,k)$ such that $f_{j,k}^\rho(z-\rho)$ is not holomorphic at $\rho$. 
Then the  family $(f_{j,k}^\rho(z-\rho))_{(j,k)\not\in\Irho}$  is a basis of the space of solutions of $(\frac{\dd}{\dd z})^N\circ\calD$ holomorphic at $\rho$. Using Eq. \eqref{eqcaractsrho} we deduce that 
  the family $(f_{j,k}^\rho(z-\rho))_{(j,k)\in\Irho}$  is a basis of $\Srho$.  If $\calD \in \K[z, \frac{\dd}{\dd z}]$,  $\rho\in\K$, and $t_j^\rho\in\K$ for any $\rho$ and any $j$, 
  then one may choose all $ g_{j,k }^\rho$ with Taylor coefficients in $\K$; moreover if $\calD$ is a $G$-operator then the  Andr\'e-Chudnovski-Katz Theorem shows that one may choose all $ g_{j,k }^\rho$ to be $G$-functions, and in this case all exponents $t_j^\rho$ are rational numbers.

Now recall  from \S \ref{subsec12} that 
$$y_{\al,i}(z) = \sum_{n=0}^\infty \frac{1}{i!} \frac{\dd^{i}}{\dd y^{i}}\Big(\frac{\Gamma(1-\{y\})}{\Gamma(-y-n)}\Big)_{|y=\al } z^n$$
and that $\hada$ denotes Hadamard's  product of formal series in $z$. With this notation we let 
$$\etanv_{j,k}^\rho(1/x) = \sum_{m=0}^k (y_{t_j^\rho,m}\hada g_{j,k-m}^\rho)(1/x) \in \C[[1/x]]$$
for any $1\leq j \leq  J(\rho)$ and $0\leq k \leq K(j,\rho)$; it is not difficult to see that $\etanv_{j,k}^\rho(1/x) = 0$ if $(j,k)\not\in\Irho$. Now 
$\oplus \Lrho : \oplus \Srho  
\stackrel{\sim}{\longrightarrow}  \oplus_{\rho } \Sbarinfrho$ is a bijective linear map (see \S \ref{secdiagLG}), so 
  that the functions $\Lrho( f_{j,k}^\rho (z-\rho))$, 
for $\rho\in\Sigma$ and  $(j,k)\in\Irho$, form a basis of formal solutions of $\calFbar \calD$ at infinity.  We claim that they can be written explicitly as follows:
 \begin{equation}\label{eqsolneta}
\Lrho( f_{j,k}^\rho (z-\rho))  = 
e^{\rho x} x^{-t_j^\rho-1}\sum_{k'=0}^k \etanv_{j,k-k'}^\rho(1/x) \frac{(\log(1/x))^{k'}}{k'!} .
\end{equation}
Moreover:
\begin{enumerate}
\item[$(i)$] If $\K$ is a subfield of $\C$ such that $\calD\in\K[z,\frac{\dd}{\dd z}]$, $\Sigma\subset\K$,  and all exponents of $\calD$ at all finite singularities  belong to $\K$, then all power series $\etanv_{j,k}^\rho(1/x) $ belong to $\K[[1/x]]$.
\item[$(ii)$] If $\calD$ is a $G$-operator then the $ g_{j,k}^\rho$ can be chosen to be $G$-functions, and then all $\etanv_{j,k}^\rho $ are \antie-functions.
\end{enumerate}

\bigskip

Indeed   Eqns. \eqref{eqdeffjk} and \eqref{eqcalclrho} yield
\begin{align*}
\Lrho( f_{j,k}^\rho (z-\rho)) 
&=\sum_{k'=0}^k  e^{\rho x}  x^{-t_j^\rho-1} \sum_{\nu=0}^{k'} (  y_{t_j^\rho,k'-\nu}\hada g_{j,k-k'}^\rho)(1/x) \frac{(\log(1/x))^{\nu}}{\nu !} \\
&=  e^{\rho x}  x^{-t_j^\rho-1} \sum_{\nu=0}^{k }  \etanv_{j,k-\nu}^\rho(1/x) \frac{(\log(1/x))^{\nu}}{\nu !}
\end{align*}
so that \eqref{eqsolneta} holds.  The rationality property $(i)$ follows at once from Proposition \ref{proprat}, and $(ii)$ is a consequence of the proof of Proposition \ref{proparith}.

\bigskip

At last, given $\rho\in\Sigma$  we have $\Card\,  \Irho = \dim\Srho = \mrho$; this is exactly the number of  functions $\etanv_{j,k}^\rho $ associated with $\rho$. Combining them as $\rho$ varies yields $\delta = \sum_\rho \mrho$  functions $\etanv_{j,k}^\rho $. The matrix $\Gamma_\infty$ of Theorem \ref{thYA} $(ii)$ is block-diagonal; with each $\rho\in\Sigma$ is associated an upper triangular block $\Gamma_{\infty,\rho}\in M_{\mrho}(\Q)$ with diagonal coefficients $t_j^\rho+1$. This completes the proof of Theorem \ref{thYA}.

\begin{Remark}
The basis \eqref{eqsolneta} is used implicitly in \S 4.2 of \cite{ateo} to compute Stokes constants.
\end{Remark}

\subsection{Arithmetic nature of the coefficients} \label{subsecarith}

Let $\calD \in \Qbar[z, \frac{\dd}{\dd z}]$ be a $G$-operator. As in \S \ref{subsecnvpreuve} we    denote by $\Sigma$ the set of all finite singularities of $\calD$, and we fix  a real number $\theta$ such that $\theta\not\equiv -  \arg(\rho-\rho') \bmod \pi$ whenever $\rho,  \rho' \in \Sigma$ are distinct. We use the same  notation as in \S \ref{secdiagLG}. The point of Theorem \ref{theocoeffsde} below is to determine the arithmetic nature of the coefficients of the matrices of all linear maps that appear in the commutative diagram of \S \ref{secdiagLG}. With this aim in view we choose ``algebraic'' bases of the vector spaces as follows.

We fix a basis of $\Srho$ for any $\rho\in\Sigma$  (resp. of $\Sinf$) consisting in elements of $\NGA\{z-\rho\}^{\Qbar}_0$ (resp. $\NGA\{1/z \}^{\Qbar}_0$), and in  the same way  a basis of  $\Sbarzero$   (resp. of $\Sbarrho$ for any $\rho\in\Sigma$) consisting in elements of $\NGA\{x\}^{\Qbar}_{-1}$ (resp. $e^{\rho x}  \NGA\{x\}^{\Qbar}_{1}$).
We also choose the basis of $\Sbarth$ in such a way that $\tay $ is represented by the identity matrix.

\begin{theo} \label{theocoeffsde} Assume   that $\calD$ is a $G$-operator. Then  in these bases:
\begin{itemize}
\item The  determinant of the matrix of $\asy$ is a   non-zero algebraic number.
\item The matrices of $\asy$, $\Vloc$, $\Vinf$,  and their inverse matrices have entries in ${\bf S}$.
\item The determinants of the matrices of $\Vloc$ and $\Vinf $ are products of values of $\Gamma$ at rational points, multiplied by non-zero algebraic numbers.
\item The matrices of $\pro$ and $\pro^{-1}$  have entries in $\GG$, and their determinants are units of $\GG$.
\end{itemize}
\end{theo}

 The first assertion follows from the results of \cite{YA1} (see the proof below). The main tool in the proof of the other ones is the commutative diagram of \S \ref{secdiagLG}, in which the maps $\Lseul$ and $\Cseul$ are explicit.  We shall also use the following restatement of Theorem 2 of \cite{gvalues}.

\begin{theo} \label{thgvalues}
Let $\calD $  be a $G$-operator of order $\mu$, and $\xi_1,\xi_2\in\Qbar\cup\{\infty\}$. Let $g\in \NGA\{z-\xi_1\}^{\Qbar}_0$ be a solution of $\calD$, and $(g_1,\ldots,g_\mu)$ be a basis of solutions of $\calD$ in $\NGA\{z-\xi_2\}^{\Qbar}_0$. Let $\varpi_1,\ldots,\varpi_\mu$ be the connection constants such that a given analytic continuation of $g(z)$ to a small cut disk \eqref{eqdisquecoupe} centered at $\xi_2$ is equal to $\varpi_1 g_1 +\ldots + \varpi_\mu g_\mu$. Then we have $\varpi_1,\ldots,\varpi_\mu\in\GG$.
\end{theo}

In this result, $z-\xi_i$ should be understood as $1/z$, and  \eqref{eqdisquecoupe} as \eqref{domaineholonilsinf}, if $\xi_i = \infty$.

\bigskip

Let us prove Theorem \ref{theocoeffsde} now.
Let $\matpro$,  $\matasy$,  $\matVinf$, $\matVloc$,$\matCinf$, $  \matCzero$  denote the respective matrices of $\pro$, $\asy$, $\Vinf $, $\Vloc$,   $\Cinf$, $\Czero$ in the bases we have fixed.  

\medskip

To begin with, let us prove the statements about $\matasy$. Let $(F_1(x), \ldots, F_\delta(x))$ and  $(H_1(x), \ldots, H_\delta(x))$ denote the bases of $\Sbarzero$ and $\oplus_{\rho\in\Sigma} \Sbarrho$ (respectively) that we have chosen, so that $F_j(x) \in \NGA\{x\}^{\Qbar}_{-1}$ and $e^{-\rho_j x } H_j(x) \in \NGA\{1/x\}^{\Qbar}_{ 1}$ for any $j$, with $\rho_1,\ldots,\rho_\delta\in\Sigma$. For any $j$, the asymptotic expansion of $F_j(x) $ as $|x|\rightarrow\infty$ in a large sector bisected by $\theta$ can be written as $\varpi_{1,j}H_1(x) + \ldots  + \varpi_{\delta,j}H_\delta(x) $; then $\matasy$ is the matrix $[\varpi_{i,j}]_{1\leq i,j \leq \delta}$. We have proved  \cite{ateo} that $\varpi_{i,j}\in{\bf S}$; we shall prove now that $\det\matasy\in\Qbar\etoile$. This implies that $\matasy^{-1} =(\det\matasy)^{-1} \tra {\rm Com } \matasy$ has coefficients in ${\bf S}$, and determinant $(\det\matasy)^{-1} \in\Qbar\etoile$.

\bigskip

Let $$
w(F_1, \ldots,F_\delta) = 
\left(\begin{matrix}
F_1 & F_1   ' &\cdots & F_1^{(\delta-1)}
\\
F _2 & F_2 ' &\cdots & F_2^{(\delta-1)}
\\
\vdots & \vdots &\cdots & \vdots
\\
F_\delta & F_\delta ' &\cdots & F_\delta^{(\delta-1)}
\end{matrix}\right)
$$ 
denote the wronskian matrix built on $(F_1 , \ldots, F_\delta )$, and $W(F_1, \ldots,F_\delta)$ be its determinant. 
Then the  asymptotic expansion of $w(F_1, \ldots,F_\delta)$ as $|x|\rightarrow\infty$ in a large sector bisected by $\theta$ is $  \matasy w(H_1, \ldots,H_\delta)$, and therefore that of $W(F_1, \ldots,F_\delta)$ is   $\det(\matasy )   W(H_1, \ldots,H_\delta)$.  Now each (formal or convergent) wronskian determinant $W(x)$ built on a basis of 
(formal or convergent) solutions of $\calFbar\calD$ satisfies   $Q_0(x) W'(x) +Q_1(x) W(x)=0$, where
$$\calFbar \calD = Q_0(x) (\frac{\dd}{\dd x})^\delta + Q_1(x) (\frac{\dd}{\dd x})^{\delta-1}  + \ldots + Q_\delta(x)$$
with $Q_0,\ldots,Q_\delta\in\Qbar[x]$ and $Q_0\neq 0$; for simplicity we assume the leading coefficient of $Q_0$ to be 1.
 Since $\calFbar \calD $ is an $E$-operator, we have $Q_0(x) = x^\mu$ where $\mu$ is the degree of  $\calFbar \calD $, and $Q_1(x) = ax^\mu + bx^{\mu-1}$ with $a,b\in\Qbar$ (see \S 5.1 of \cite{YA1}).  Therefore $W'(x) + (a+\frac{b}{x} )W(x) = 0$ so that $W(x) = c_W x^b e^{ax}$ where $c_W\in\C\etoile$ depends on $W$. Applying this property to both wronskians introduced previously, we obtain that $\det(\matasy )   c_H x^b e^{ax}$ is the  asymptotic expansion of  $  c_F x^b e^{ax}$   as $|x|\rightarrow\infty$ in a large sector bisected by $\theta$; this means  $ c_F = \det(\matasy )   c_H $. Now we have $W(F_1, \ldots,F_\delta)\in  \NGA\{x\}^{\Qbar}_{-\delta}$ and $e^{-\rho_1 x -\ldots-\rho_\delta x } W(H_1, \ldots,H_\delta)\in  \NGA\{1/x\}^{\Qbar}_{ \delta}$  so that $\rho_1+\ldots+\rho_\delta = a$ and $c_F,c_H\in\Qbar\etoile$. Therefore  $\det(\matasy ) = c_F /  c_H\in\Qbar\etoile$.

\medskip

Now since $\GG$ and ${\bf S}$ are $\Q$-vector spaces, the truth of Theorem \ref{theocoeffsde} is independent from the bases chosen, as long as the $\Qbar$-structures are preserved (see for instance \S 8 of \cite{Boualg}).  Therefore we may assume that the basis of $\oplus_\rho \Srho$ is $( f_{j,k}^\rho(z-\rho) )_{\rho\in\Sigma, (j,k)\in\Irho}$ with the notation of \S\ref{subsecnvpreuve}, ordered   lexicographically with respect to $(\rho,j, k)$. 
Let   us consider  an element of the basis of $\Sinf$ we have chosen. Using Eq. \eqref{eqlemdim} proved in \S \ref{subsec32}
 and the Andr\'e-Chudnovski-Katz Theorem,  
  it is represented by a function $f$ holomorphic on the cut plane $\cut$,   such that  $(\frac{\dd}{\dd z})^N  \calD f  = 0$ and $f\in \NGA\{1/z\}^{\Qbar}_{ 0}$ around $\infty$.
  Recall from \S \ref{subsecnvpreuve} that for any $\rho\in\Sigma$, the functions $ f_{j,k}^\rho(z-\rho)$ with $1\leq j \leq J(\rho)$ and $0\leq k \leq K(\rho,j)$ form a basis of solutions of $(\frac{\dd}{\dd z})^N  \circ \calD$ in $ \NGA\{z-\rho\}^{\Qbar}_{ 0}$. Since $(\frac{\dd}{\dd z})^N  \circ \calD$  is a $G$-operator, we may expand $f$ in this local basis at $\rho$, with connection constants $\varpi_{j,k}^\rho$ in $\GG$ (using Theorem \ref{thgvalues}). Omitting the indices $(j,k)\not\in\Irho$ and letting $\rho$ vary, we obtain the coefficients in the column of $\matpro$ corresponding to $f$. This concludes the proof that $\matpro$ has coefficients in $\GG$.

\medskip

Let us focus on $ \matCinf$ now. Since the basis $( f_{j,k}^\rho(z-\rho) )_{\rho\in\Sigma, (j,k)\in\Irho}$ is ordered   lexicographically with respect to $(\rho,j, k)$,   the matrix $ \matCinf$ is upper triangular and its diagonal coefficients are $\Gh(1-\{t_j^\rho\})$, where $t_j^\rho \in \Q$ appears in Eq. \eqref{eqdeffjk}; therefore
$$\det \matCinf = \prod_{\rho\in\Sigma} \prod_{j=1}^{J(\rho)}   
\Big(\Gh(1-\{t_j^\rho\})\Big)^{K(\rho,j)+1}$$
where $\{t_j^\rho\}$ is the fractional part of $t_j^\rho$ and $\Gh(s) = 1/\Gamma(s)$.
In this product some factors should have been omitted, namely those which correspond to triples $(\rho,j, k)$ for which $f_{j,k}^\rho(z-\rho) $ is holomorphic at $\rho$; however in this case $\Gh(1-\{t_j^\rho\}) = 1$ so  that equality holds anyway.

The situation is the same if $\rho$ is replaced with $\infty$, upon replacing $z-\rho$ with $1/z$. The only difference is that $(1/z)^{t_j^\infty}$ should be seen as $z^{-t_j^\infty}$ when applying $\Linf$, so that we obtain
$$\det \matCzero =  \prod_{j=1}^{J(\infty)}   
\Big(\Gh(1-\{-t_j^\infty\})\Big)^{K(\infty,j)+1}.$$

Now the  matrices of $\oplus \Lrho$ and $\Linf$ have algebraic coefficients (see Proposition \ref{proprat}), and we have proved that $\det\matasy\in\Qbar\etoile$, so that 
$$\det \matVloc  = \frac{\cloc}{\det \matCinf } = \cloc  \prod_{\rho\in\Sigma} \prod_{j=1}^{J(\rho)}   
\Big(\Gamma(1-\{t_j^\rho\})\Big)^{K(\rho,j)+1} $$
and
$$\det \matVinf  = \frac{\cinf}{\det \matCzero} = \cinf  \prod_{j=1}^{J(\infty)}   
\Big(\Gamma(1-\{-t_j^\infty\})\Big)^{K(\infty,j)+1} $$
with $\cloc,\cinf\in\Qbar\etoile$. Moreover $ \matVloc  =   \matpro \matVinf$ (see \S \ref{subseclaplace}) so that
 $\det\matpro = \frac{\cloc}{\cinf} \frac{\det \matCzero}{\det \matCinf }$. Now we observe 
  the following fact:   if $x_1,\ldots,x_p,y_1,\ldots,y_q\in\Rplusetoile$ satisfy  $x_1+\ldots+x_p = y_1+\ldots+y_q$ then
$$\frac{\prod_{i=1}^p \Gamma(x_i)}{\prod_{j=1}^q \Gamma(y_j)} = \frac{\prod_{i=1}^{p-1} B(x_1+\ldots+x_i, x_{i+1})}{\prod_{j=1}^{q-1} B(y_1+\ldots+y_j, y_{j+1})} $$
since $B(a,b) = \frac{\Gamma(a)\Gamma(b)}{\Gamma(a+b)}$. Since $\calD$ is Fuchsian, 
  Fuchs' relation on exponents yields 
$$   \sum_{\rho\in\Sigma\cup\{\infty\}} \sum_{j=1}^{J(\rho)}   (K(\rho,j)+1) t_j^\rho \in \Z $$
so that using this fact we obtain that the number $\frac{\det \matCinf}{\det \matCzero}$ is a quotient of products of values $B(x,y)$ with $x,y\in \Qplusetoile$;  in particular it is a unit of $\GG$ (see Proposition 1 of \cite{gvalues}). Therefore $\det\matpro = \frac{\cloc}{\cinf} \frac{\det \matCzero}{\det \matCinf }$  is also a unit of $\GG$, and   $\matpro^{-1} =(\det\matpro)^{-1} \tra {\rm Com } \matpro$ has coefficients in $\GG$.

This concludes the proof of  Theorem \ref{theocoeffsde}.

\subsection{Special values and exponential periods} \label{subsec43}

In this section we prove the arithmetic results stated in the introduction and explain the connection with exponential periods.  To deduce Theorem \ref{thnvcor} from Theorem \ref{theocoeffsde} proved in \S \ref{subsecarith}, we just have to choose an $E$-operator that annihilates $\ae(1/x)$, using Theorem 4.6 of \cite{YA1}. We obtain  also Corollary \ref{corintrode} at once, since any Stokes matrix is the matrix of 
$ {\asy}\circ  {\asypr}^{-1}$ in the basis $\Bbarinf$ (with the notation of the introduction).

Before we prove Theorem \ref{thintrotr}, let us study briefly the set \ensembleae of all numbers 
$\ae_\theta(\xi)$ where $\xi\in\Qbar\etoile$  and $\ae$ is an \antie-function; here 
$\theta = \arg\xi$ and  $\ae_\theta =  \asymu \ae$ is Ramis' 1-summation of $\ae$ in the direction $\theta$ if $\theta$ is not anti-Stokes, and  $\ae_\theta = \asyepsmu\ae$ for any small $\eps>0$ if $\theta$ is  anti-Stokes (this is independent from the choice of such an $\eps$). 
 From now on, we shall always restrict to the case of non-anti-Stokes directions, because the opposite case can be proved along the same lines.  We shall also assume, for simplicity, that $\ae(1/x) = \sum_{n=1}^\infty a_n x^{-n} $ has no constant term. Then its  formal Borel transform $g(z) = \sum_{n=1}^\infty \frac{a_n}{(n-1)!} z^{n-1} $   is a $G$-function, and we have  
\begin{equation} \label{equnsommmetath}
 \asymu  \ae(1/x) = \int_0 ^{e^{i \theta} \infty} g(z) e^{-xz} \dd z
\end{equation}
provided $|\arg(1/x) -  \theta| < \frac{\pi}{2}$ (see \cite{Ramis}), in particular for $1/x = \xi$. Moreover letting $\tilde{\ae}(x) = \ae(\xi x)$ we have $\tilde{\ae}_0(1) = \ae_\theta(\xi)$ so that we may restrict to $\xi=1$  in the definition of \ensembleaesansesp; therefore \ensembleae is a ring.
 Moreover any $G$-function $g(z)$ is the  formal Borel transform of an  \antie-function with  no constant term 
so that all numbers \eqref{equnsommmetath} belong to \ensembleaesansesp, provided $x$ is algebraic and  $\theta = \arg(1/x)$ is not anti-Stokes. In particular  \ensembleae contains Gompertz' constant $\int_0 ^{+\infty}  \frac{e^{-t}}{1+t} \dd t$, and $\sqrt{\pi} {\rm Ai}(z)$ for any $z\in\Qbar$ where ${\rm Ai}(z)$ is Airy's oscillating integral (see \cite{YA1}).

\bigskip

Let us come now to the properties of $\calV$. Let $\xi\in\Qbar$, $s\in\Q$, and $h\in\Qbar\{z\}^A_s$ be a   Gevrey series of order $s $ of arithmetical type. If $s=0$ then $h$ is a $G$-function so that $h(\xi) \in\GG\subset{\bf S}\subset\calV$. If $s\neq 0$ then $h(z^{|s|})\in\NGA\{z\}_{\eps(s)}$, where $\eps(s) \in\{-1,1\}$ is the sign of $s$ (see Proposition 1.4.1 of \cite{YA1}). Using Theorem \ref{thnvcor} if $s>0$, we obtain also $h(\xi)\in\calV$. This proves that $\calV$  contains all values at algebraic points of 
  Nilsson-Gevrey series  of any order   of arithmetical type with algebraic coefficients $\lambda_{\al,j,k}$ (up to 1-summation in any direction in the case of divergent series), and yields  \ensembleae $\subset\calV$. To deduce  Theorem \ref{thintrotr}, we denote by $\calV'$ the ${\bf S}$-module generated by the numbers $e^\rho \chi$, with $\rho\in\Qbar$ and $\chi\in$ \ensembleaesansesp. Since \ensembleae $\subset\calV$ and $e^\rho\in\E\subset\calV$ for any $\rho\in\Qbar$, we have  $\calV'\subset\calV$. On the other hand, given $\xi\in\Qbar$ and an  $E$-function $F(x)$ there exist $A\geq 1$, $\varpi_1,\ldots,\varpi_A\in{\bf S}$, $\rho_1,\ldots,\rho_A\in\Qbar$ and $\ae_1(1/x),\ldots, \ae_A(1/x)\in\NGA\{1/x\}^{\Qbar}_1$ such that $\sum_{a=1}^A \varpi_a e^{\rho_a x}\ae_a(1/x)$ is the asymptotic expansion of $F(x)$ as $|x|\rightarrow\infty$  in a large sector bisected by $\theta = \arg\xi$ if it is not anti-Stokes, by $\theta+\eps$ otherwise  (see \cite{ateo}). Applying 1-summation in this direction and evaluating at  $\xi$, we obtain $\E \subset\calV'$; this concludes the proof of Theorem \ref{thintrotr}.
  
  \bigskip
  
  \begin{Remark} In the definition of \ensembleae we could have considered $\asymoinsepsmu \ae$ instead of $ \asyepsmu\ae$ if $\theta$ is  anti-Stokes. This defines also a ring \ensembleaesansesp', and all properties of \ensembleae stated in this paper hold also with \ensembleaesansesp'. However we have not been able to prove that \ensembleaesansesp' = \ensembleaesansesp.
  \end{Remark}
  
  \bigskip
  
  To conclude this discussion, we observe that $\Gamma(r) = \int_0^{ \infty} t^{r-1}e^{-t} \dd t$ (with $r\in\Q\setminus\Zneg$) and Euler's constant $\gamma = -\int_0^\infty \log(t)e^{-t}\dd t$ are exponential periods. Therefore the conjectural inclusion $\GG\subset\calP[1/\pi]$ implies ${\bf S} \subset \calP_e[1/\pi]$. Now Bombieri-Dwork conjecture predicts that $G$-functions come from geometry, so that Eq.  \eqref{equnsommmetath} and Theorem \ref{thintrotr} suggest the inclusion $\calV\subset \calP_e[1/\pi]$. In view of Andr\'e's results (see \cite{Andre} and \cite{AndreLMS}), it seems reasonable to believe that 
$\calV$ could be equal to $\calP_e[1/\pi]$. 

\subsection{An example related to Gompertz' s constant} \label{subsec44}

In this section we work out a specific example involving Gompertz' s constant $\int_0 ^{+\infty}  \frac{e^{-t}}{1+t} \dd t$. The notation is not always recalled, but it is the same as in the rest of the paper. Unless otherwise stated, all properties can be proved as special cases of results given in the present paper. Other references on this example include \cite{LodayGazette} and Proposition 1 of \cite{Michigan}.

\medskip

Let $\calD = z(1-z) \frac{\dd}{\dd z}-z$ and $g_1(z) = \frac1{1-z}$. Then $\calD\in\Qbar[z,\frac{\dd}{\dd z}]$ is of minimal order $\mu=1$ among the differential operators such that $\calD g_1=0$. Since $g_1$ is rational, it is a $G$-function and $\calD$ is a $G$-operator. The set $\Sigma$ of finite singularities of $\calD$ is $\{0,1\}$. We fix $\theta$ such that $-\pi < \theta< \pi$ and $\theta\neq 0$, and work in the simply connected open set $\cut$ defined by $-\theta-\pi < \arg z < -\theta+\pi$ and  $-\theta-\pi < \arg ( z-1)  < -\theta+\pi$. Each finite singularity has multiplicity 1, so that the associated spaces ${\mathcal S}_0$ and ${\mathcal S}_1$ of microsolutions are one-dimensional. We choose $\bar g_0$ and $\bar g_1$ as respective bases of ${\mathcal S}_0$ and ${\mathcal S}_1$, where the overline denotes the image  in a quotient group and $g_0(z) = \frac{\log z}{1-z}$; we observe that $\calD g_0 =1\in \cz_{<N}$, with $N=1$ in this example. Choosing $(\bar g_0,\bar g_1)$ as a basis of $\Sinf$, the matrix of $\pro$ is the identity matrix.

\medskip

The associated $E$-operator is $\calFbar\calD = x(\frac{\dd}{\dd x})^2 + (1-x) \frac{\dd}{\dd x} - 1$. Letting $E(x) =\sum_{n=1}^\infty\frac{(-1)^nx^n}{n\cdot n!}$, a basis of $\Sbarth$ is given by the functions $e^x (E(x)  + \log x)$ and $e^x$. Both have generalized Taylor expansions at 0 in $\NGA\{x\}_{-1}^{\Qbar}$; indeed $e^xE(x)$ is an $E$-function. The  Laplace transforms of $g_0$ and $g_1$ can be written as follows (since $\pro$ is given by the identity matrix):
$${\mathcal L}_0^{-1}\Big(  \frac{\log z}{1-z} \Big) = {\mathcal L}_\infty^{-1}\Big(  \frac{\log z}{1-z} \Big) = \int_0 ^{+\infty}  \frac{e^{-t}}{x+t} \dd t = -e^x(E(x)+\log(x)+\gamma)  ,$$
$${\mathcal L}_1^{-1}\Big(  \frac1{1-z} \Big) ={\mathcal L}_\infty^{-1}\Big(  \frac1{1-z} \Big) =  e^x.$$
Letting $y(x) = -e^x(E(x)+\log(x)+\gamma) $, we obtain another basis $(y(x),e^x)$ of $\Sbarth$, with $y(x) \in \overline{{\mathcal S}}_{[\theta]}^0$ and $e^x  \in \overline{{\mathcal S}}_{[\theta]}^1$. Moreover $y(x)$ admits the asymptotic expansion $\widehat y(x) = \sum_{n=0}^\infty (-1)^n n! x^{-n-1}$ as $|x| \to\infty$ in a large sector bisected by $\theta$,  so that $\widehat y(x) $ is a basis of $\overline{{\mathcal S}}_{\infty}^0$ contained in $\NGA\{x\}_1^{\Qbar}$. On the other hand,  $e^x$  is a basis of $\overline{{\mathcal S}}_{\infty}^1$ contained in $e^x \NGA\{x\}_1^{\Qbar}$.  In the bases we have chosen, the matrix of $\asy$ studied in Theorem \ref{theocoeffsde} is $\left(\begin{matrix} -1&0\\ -\gamma&1\end{matrix}\right)$. Gompertz' s constant appears as the value at $x=1$ of the \antie-function $\widehat y(x)$. 

\begin{Remark} \label{remapp}
We have $\calD = z \calD'$ where   $\calD' =  (1-z) \frac{\dd}{\dd z}-1$ is also a $G$-operator. Eventhough $\calD $ and $ \calD'$ have the same solutions, they do not have the same microsolutions at 0. Applying our results to $\calD'$ does not yield anything related to Gompertz' constant.
\end{Remark}

\providecommand{\bysame}{\leavevmode ---\ }
\providecommand{\og}{``}
\providecommand{\fg}{''}
\providecommand{\smfandname}{\&}
\providecommand{\smfedsname}{eds.}
\providecommand{\smfedname}{ed.}
\providecommand{\smfmastersthesisname}{M\'emoire}
\providecommand{\smfphdthesisname}{Th\`ese}

\noindent S. Fischler, \'Equipe d'Arithm\'etique et de G\'eom\'etrie Alg\'ebrique, 
Universit\'e  Paris-Sud, B\^atiment 425,
91405 Orsay Cedex, France

\medskip

\noindent T. Rivoal, Institut Fourier, CNRS et Universit\'e Grenoble 1, 
100 rue des maths, BP 74, 38402 St Martin d'H\`eres Cedex, France 

\bigskip

\noindent Keywords: $E$-operator,  $G$-operator, Microsolution,  Fuchsian operator, Laplace transform, Special value, Arithmetic Gevrey series.

\medskip

\noindent Math. Subject Classification (2010): 11J91 (Primary);  33E30, 34M40,  44A10  (Secondary).

\end{document}